\documentclass[10pt,reqno]{amsart}
\usepackage{amssymb}
\usepackage{amsmath}
\usepackage{graphicx}
\theoremstyle{plain}
\newtheorem{Thm}{Theorem}
\newtheorem{Cor}[Thm]{Corollary}

\newtheorem{Mres}{Main Result}
\newtheorem{\theMain}{}
\newtheorem{Lem}[Thm]{Lemma}

\newtheorem{Prop}[Thm]{Proposition}
\theoremstyle{definition}
\newtheorem{Def}{Definition}
\theoremstyle{remark}

\def\determinant#1{\left|#1\right|}

\newcommand{\ntt}{\normalfont\ttfamily}
\newcommand{\pkg}[1]{{\protect\ntt#1}}
\DeclareMathOperator{\sgn }{sgn }

\DeclareMathOperator{\len }{len }

\DeclareMathOperator{\kr }{kr}
\begin{document}
 
\begin{titlepage}

\begin{center}
\vspace*{1.325in}
EDGE EFFECTS ON LOCAL STATISTICS IN LATTICE DIMERS:\\
A STUDY OF THE AZTEC DIAMOND (FINITE CASE)\\
\vspace{1.365626in}
BY\\

\vspace{.2in}
HARALD HELFGOTT\\
\vspace{1.11875in}
THESIS\\
\vspace{.2in} 
Submitted in partial fulfillment of the requirements\\
for the degree of Bachelor in Arts in Mathematics\\
at Brandeis University, 1998\\
\vspace{1.19750in}
Waltham, Massachusetts
\end{center}
\end{titlepage}
\address{Mathematics Department \\
  Brandeis University
  Waltham, MA 02254-9110}
\email{hhelf@cs.brandeis.edu}
\section{Introduction}

\begin{figure}
\centering \includegraphics{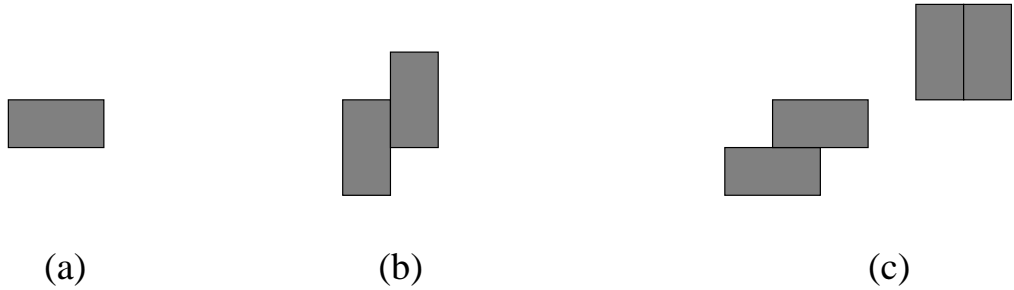}
\caption{A few examples of local patterns} \label{fig:ex}
\end{figure}

A {\em tiling} of a checkerboard with dominoes is a way of putting dominoes
on the board so that no square of the board is uncovered and no two dominoes
overlap.  Given a local pattern (see figure \ref{fig:ex} for examples), a 
location in the board, and the shape and size of the board, how many tilings
of the board have the given pattern at the given location? (Alternatively,
we can substitute ``bond'' for ``domino'' and ``particle'' for 
``square'', and ask for the probability of local patterns in a system
of particles each of which bonds with exactly one of its neighbors.)

\begin{figure}
        \begin{minipage}[b]{0.5\linewidth}
                \centering \includegraphics[height=2in]{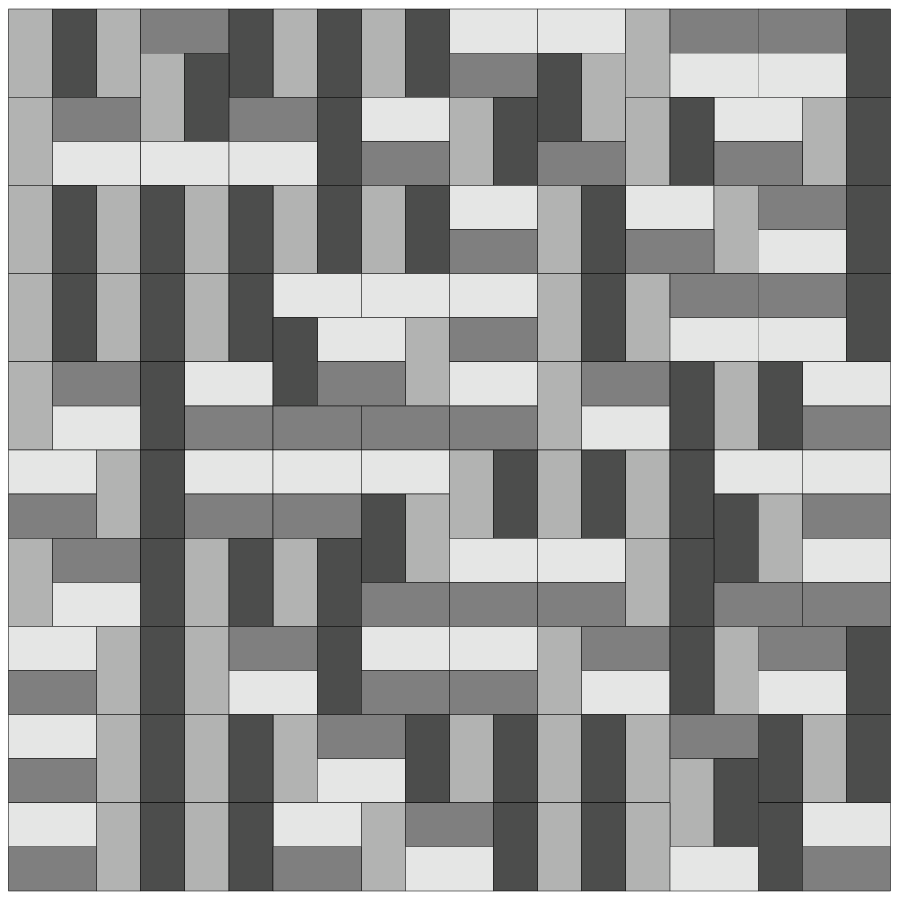}
                \caption{Random tiling of square of side 20}\label{fig:g}   
        \end{minipage}%
        \begin{minipage}[b]{0.5\linewidth}
                \centering \includegraphics[height=2in]{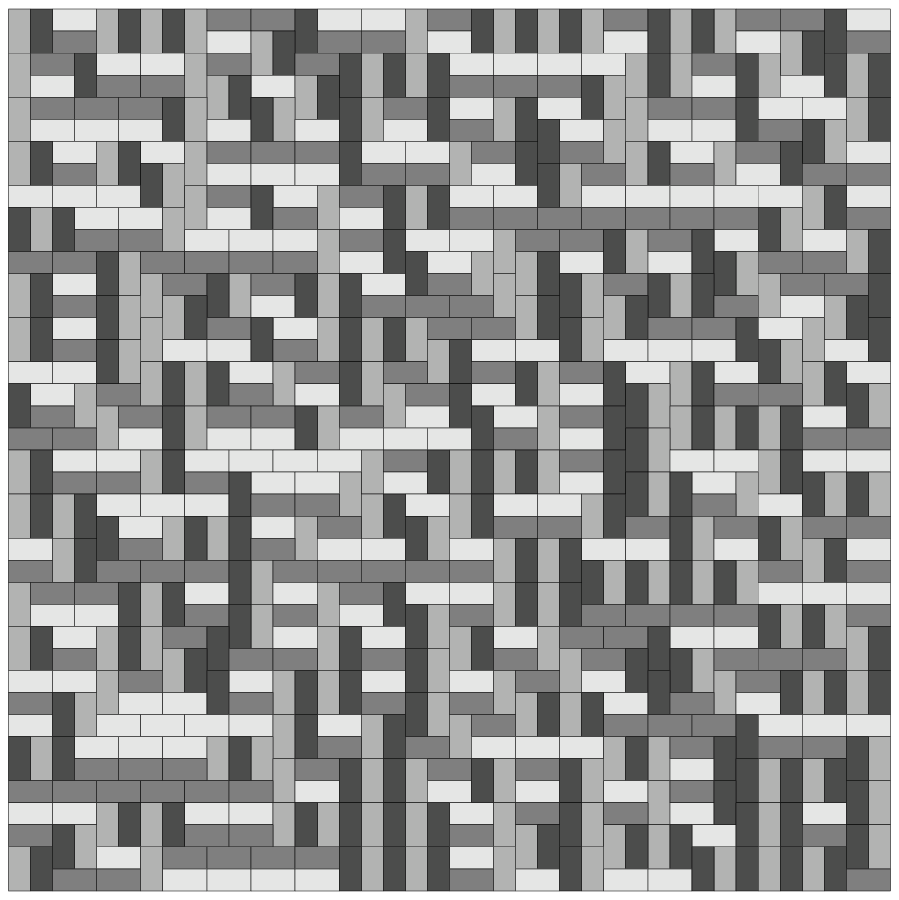}
                \caption{Random tiling of square of side 40}\label{fig:h}
        \end{minipage}%
\end{figure}

\begin{figure}
        \begin{minipage}[b]{0.5\linewidth}
                \centering \includegraphics[height=2.5in]{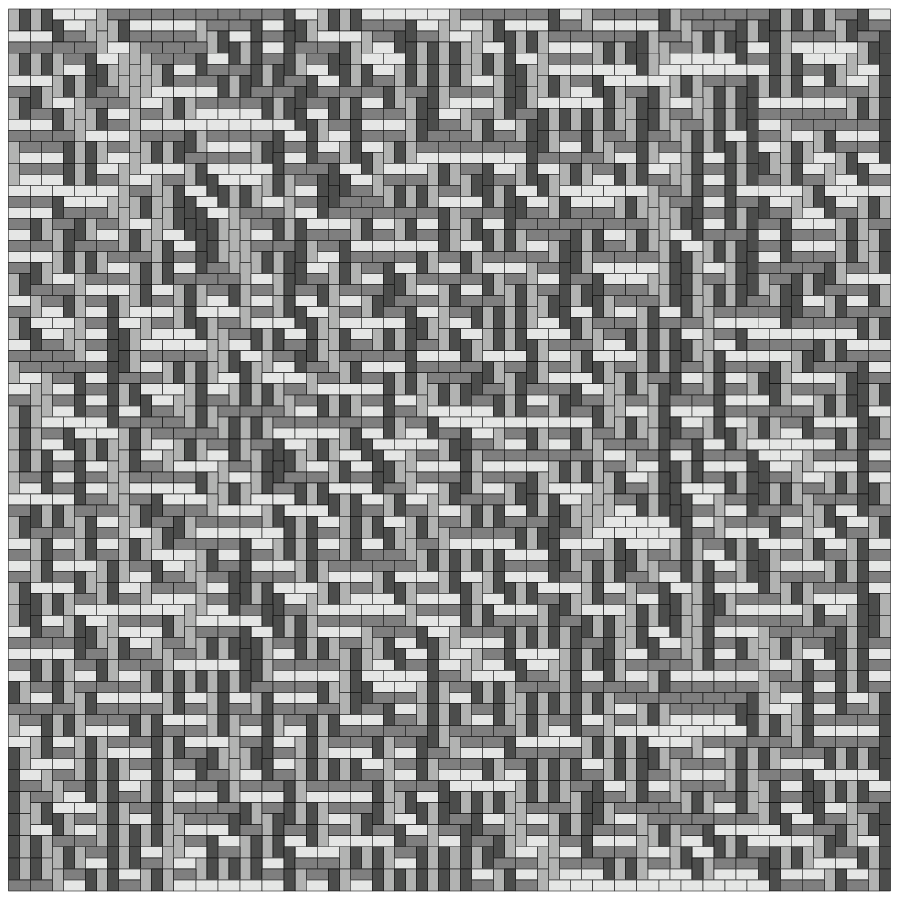}
                \caption{Random tiling of square of side 80}\label{fig:i}   
        \end{minipage}%
\end{figure}

\begin{figure}
\centering \includegraphics[height=2.5in]{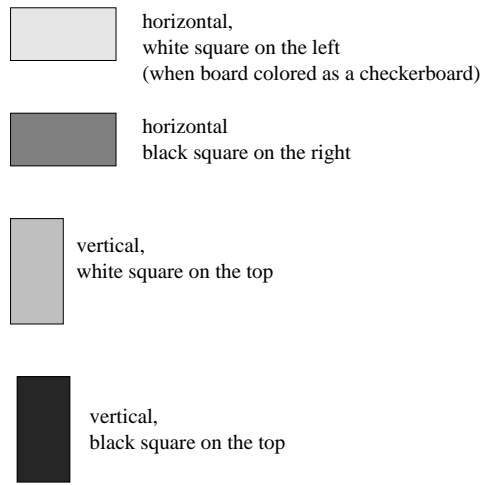} 
\caption{Shading chart} \label{fig:chart}
\end{figure}

\begin{figure}
	\begin{minipage}[b]{0.5\linewidth}
		\centering \includegraphics[height=2in]{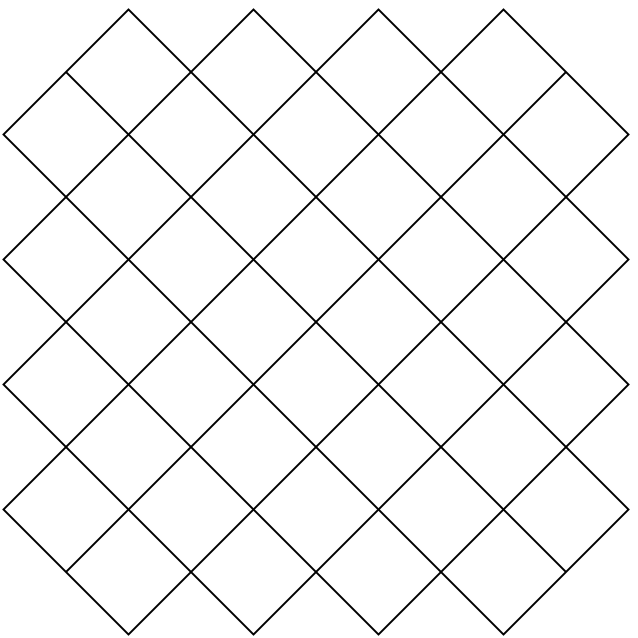}
	\caption{Aztec diamond of order $4$, as a board} \label{fig:azb}
	\end{minipage}
	\begin{minipage}[b]{0.5\linewidth}
		\centering \includegraphics[height=2in]{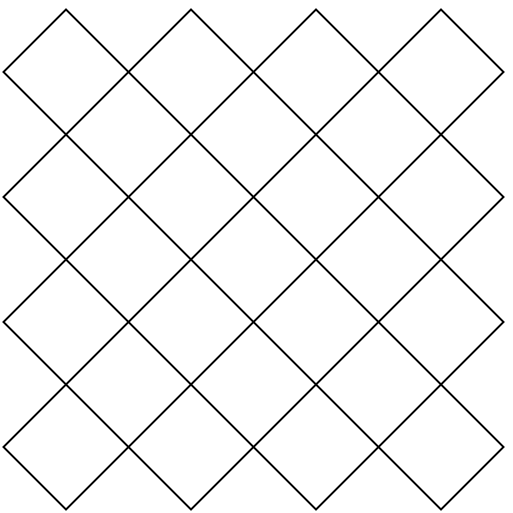}
	\caption{Aztec diamond of order $4$, as a graph} \label{fig:azg}
	\end{minipage}
\end{figure}

\begin{figure}
        \begin{minipage}[b]{0.5\linewidth}
                \centering \includegraphics[height=2in]{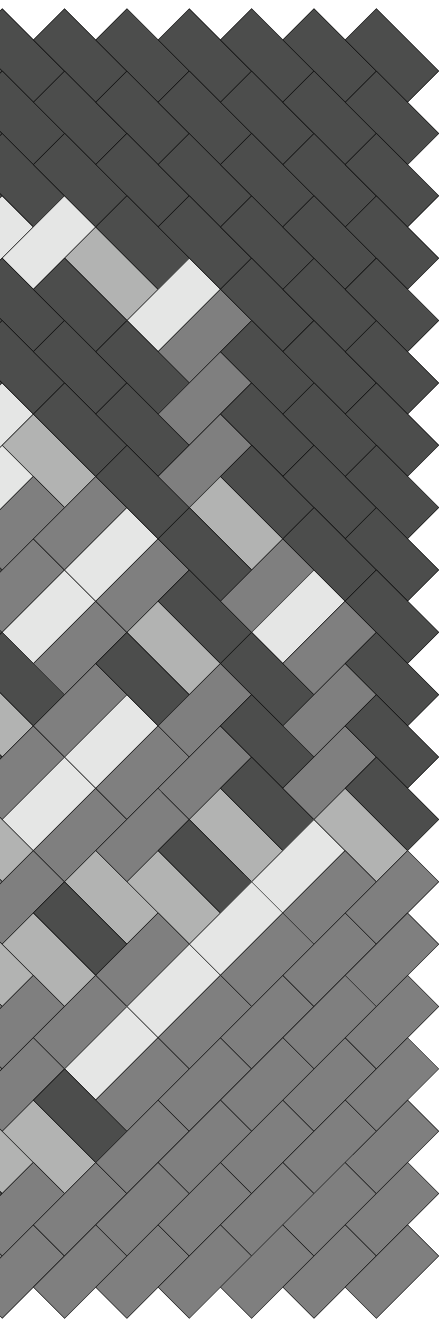}
                \caption{Random tiling of Aztec diamond of order (side) 20}\label{fig:b}   
        \end{minipage}%
        \begin{minipage}[b]{0.5\linewidth}
                \centering \includegraphics[height=2in]{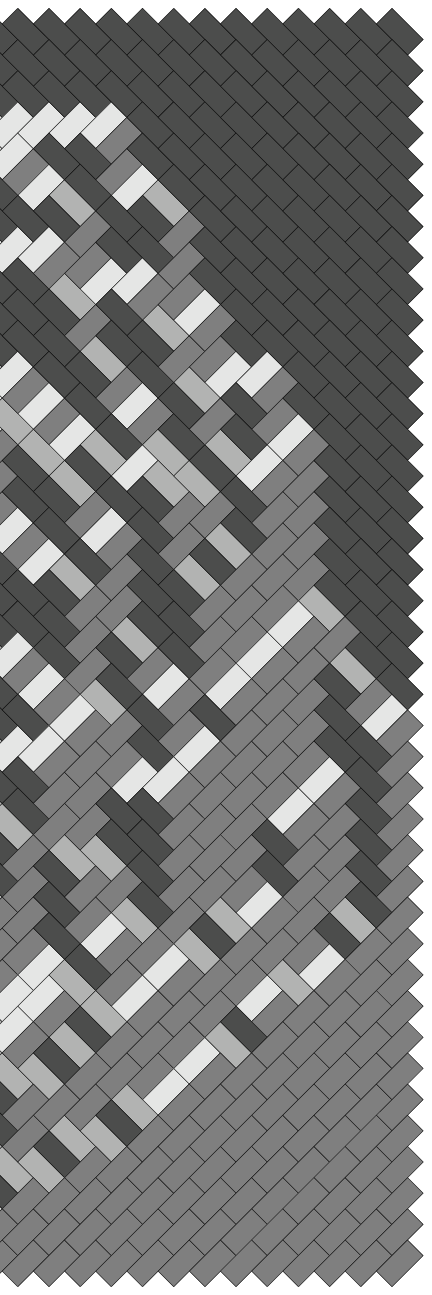}
                \caption{Random tiling of Aztec diamond of order 40}\label{fig:c}
        \end{minipage}%
\end{figure}

\begin{figure}
        \begin{minipage}[b]{0.5\linewidth}
                \centering \includegraphics[height=2in]{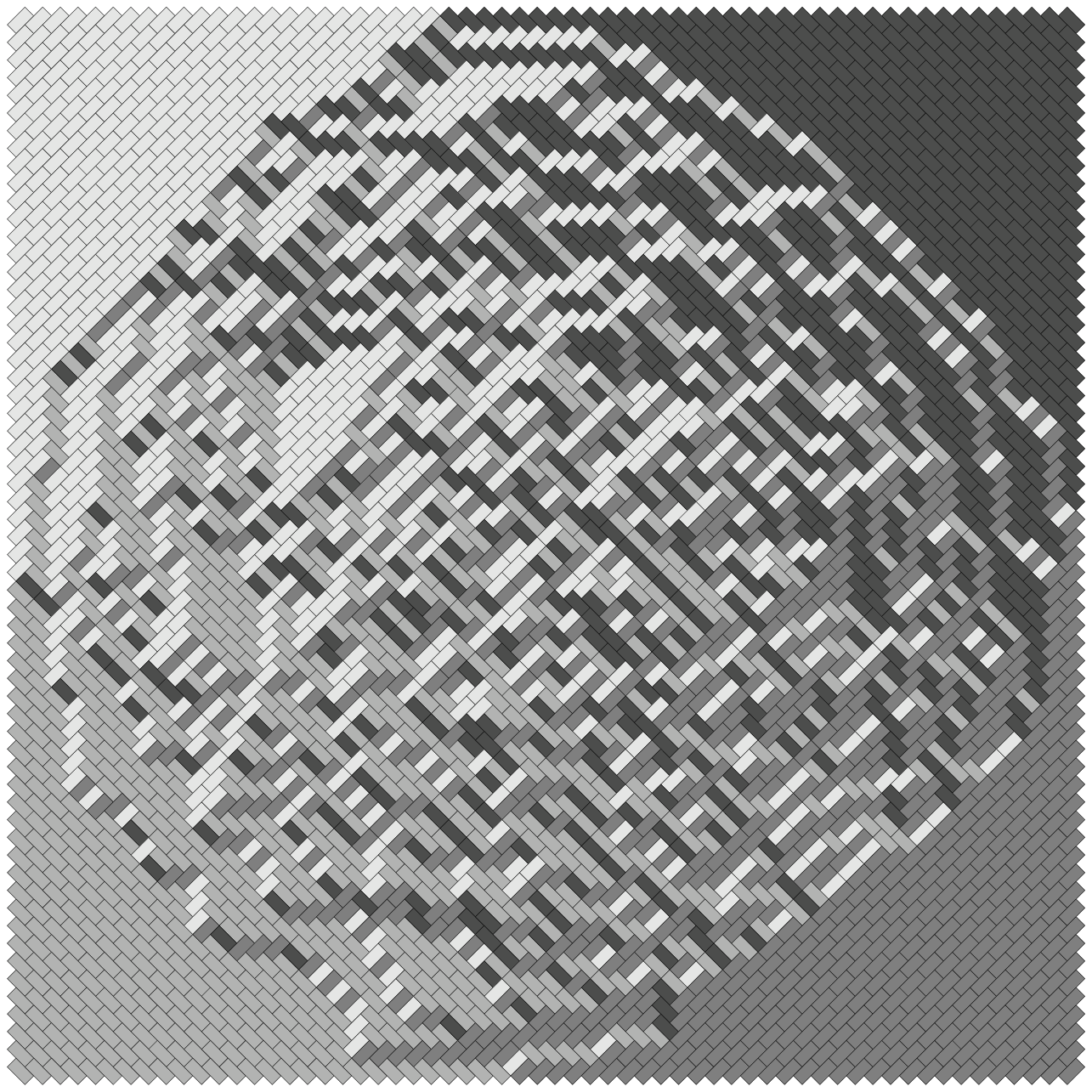}
                \caption{Random tiling of Aztec diamond of order 60} \label{fig:d}   
        \end{minipage}%
        \begin{minipage}[b]{0.5\linewidth}
                \centering \includegraphics[height=2in]{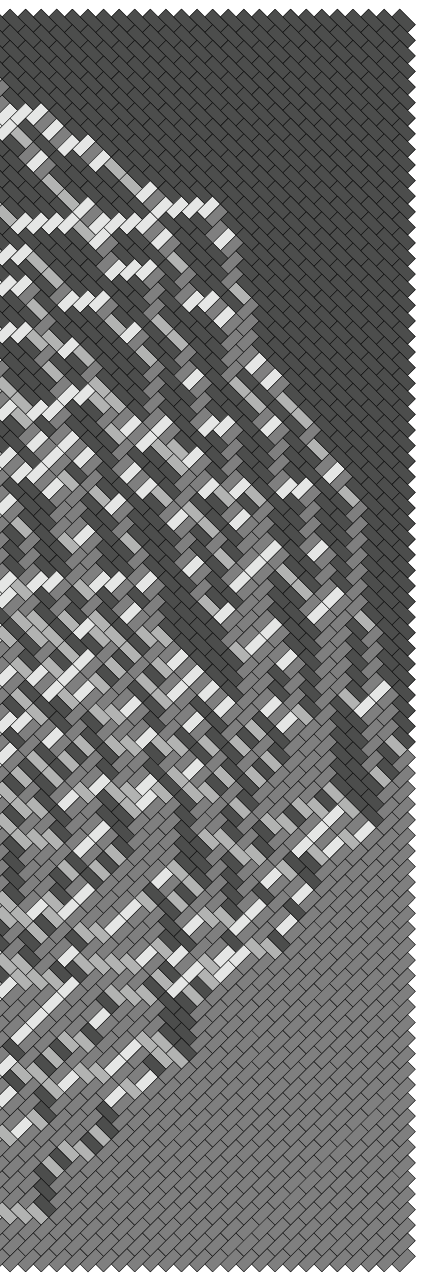}
                \caption{Random tiling of Aztec diamond of order 80}\label{fig:e}
        \end{minipage}%
\end{figure}

\begin{figure}
	\centering \includegraphics[height=4in]{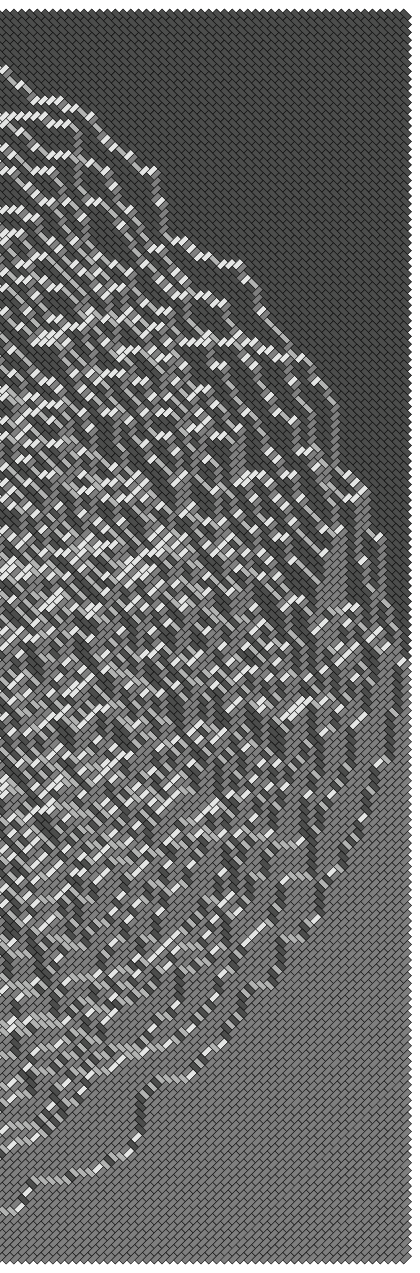}
	\caption{Random tiling of Aztec diamond of order 160}\label{fig:f}
\end{figure}

Suppose that the squares of the board are very small compared to the
board itself. For some board shapes, the probability of finding a pattern
at a given location will be the same for almost all locations. This is
the case for the square board. (See figures \ref{fig:g} to \ref{fig:i},
where tiles are colored according to their direction and parity for
the sake of clarity; see figure \ref{fig:chart} for the coloring scheme.)
There are some boards, however, for which the probability does depend
on the location. Consider, for example, the {\em Aztec diamond}, that is,
the board whose boundary is a square tilted $45$ degrees 
(figures \ref{fig:azb} and \ref{fig:azg}). 
In random tilings of the Aztec diamond, we usually find brick-wall patterns 
outside the inscribed circles, and more complicated behavior inside
the circle. (See figures \ref{fig:b} to \ref{fig:f}.)

The probabilities of local patterns in a rectangular board were computed
recently \cite{Ken}. Until now, there was no other board for which
the probabilities of all local patterns were known. Many experiments 
and some important partial results \cite{CEP} had shown that, as already
stated, the probabilities of patterns in the Aztec diamond depend on location.
This qualitative difference between the Aztec diamond and the rectangular
board made the former as worthy of analysis as the latter. 
The main result of this work
is an expression for the 
probability of any local pattern in a random tiling of the Aztec diamond.
The expression is a determinant of size proportional to the number of
squares in the pattern, just like Kenyon's expression \cite{Ken} for
the probabilities in the rectangural board,
 
\begin{Mres}
The probability of a pattern covering white squares $v_1,v_2,\dotsb v_k$
and black squares $w_1,w_2,\dotsb w_k$ of
 an Aztec diamond of order $n$ is equal
to the absolute value of
\[\determinant{c(v_i,w_j)}_{i,j=1,2,\dotsb k}.\]
The {\em coupling function} $c(v,w)$ at white square $v$ and black
square $w$ is 
\[2^{-n} \sum_{j=0}^{x_i-1} \kr(j,n,y_i-1) 
			    \kr({y\prime }_i - 1,n-1,n-(j+{x\prime }_i-x_i))
\]
for ${x\prime }_i > x_i$ and
\[-2^{-n} \sum_{j=x_i}^n \kr(j,n,y_i-1) 
			 \kr({y\prime }_i-1,n-1,n-(j+{x\prime }_i-x_i))
\]
for ${x\prime }_i \leq x_i$, where $(x_i,y_i)$ and $({x\prime}_i,{y\prime}_i)$
are the coordinates of $v$ and $w$, respectively, in the coordinate
system in figure \ref{fig:coor}, and the {\em Krawtchouk polynomial}
$\kr(a,b,c)$ is the coefficient of $x^a$ in $(1-x)^c\cdot (1+x)^{b-c}$.
\end{Mres}

Our line of attack is as follows.
\begin{enumerate}
\item Reduce the problem of finding probabilities of patterns to an 
enumerative problem;
\item Reduce the enumerative problem to a simpler one involving
Aztec diamonds with two holes rather than arbitrary even-area holes;
\item Compute the weighted number of tilings of an Aztec diamond with two holes.
\end{enumerate}

The first two steps involve known techniques, and were already considered
to be a plausible strategy by other researchers. The third step is new.

Henry Cohn is currently analyzing the case of the board with infinitely small
squares by approximating the sum of Krawtchouk polynomials in our main
result as an integral for $n\to \infty$. 
His results will be presented in a later, joint version
of this paper. 

\begin{figure}
	\centering \includegraphics[height=1in]{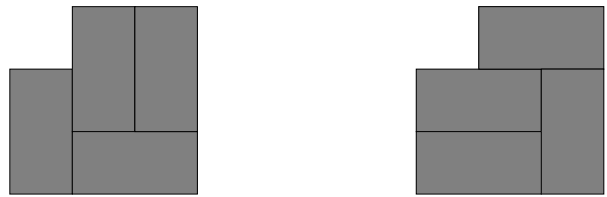}
	\caption{These two patterns
	have the same probability of being found in a random pattern
at any given place} \label{fig:nomatter}
\end{figure}



\section{The Kasteleyn Matrix}

What is the probability of finding a pattern in a tiling of a given board?
It is equal, by definition, to the number of tilings of the
board with the given pattern, divided by the total number of tilings
of the board. Clearly, the number of tilings of the board with the
given pattern depends only on the squares occupied by the pattern, and not
by how they are tiled by the pattern: there is a one-to-one correspondence
between tilings with the given pattern, and tilings of the
board with the squares covered by the pattern removed. (See figure 
\ref{fig:nomatter}.) Thus, what we want to know is the number of tilings
of the board with the squares covered by a pattern removed, divided
by the number of tilings of the board.

Kasteleyn \cite{Kast} showed that the number of tilings of a board can
be expressed as the absolute value of a determinant with half as many rows
as there are squares in the board. Thus, in a sense, our problem is already
solved: since determinants can be computed in time polynomial on the
number of rows, we can compute the probability of any pattern in a board
in time polynomial on the size of the board. Unfortunately, there are two
disadvantages to this approach. The first one is that, if the Aztec
diamond in question has side $n$, we will have to compute a determinant
of side $\frac{n^2}{2}$, and that demands a considerable amount
of time and space. Even more seriously, a
sequence of determinants of varying size is very hard to analyze 
asymptotically. We would like to analyze such a sequence in order to
determine the probability of a pattern near a point for Aztec diamonds
of very high order, that is, of very fine ``grain'' (see figure
\ref{fig:f}).
Thus, a Kasteleyn determinant is not good enough.
A determinant whose size depended only on the size of the pattern,
and not on the size of the board, would be much easier to manipulate.

In this section, we will prove that the number of tilings of a board
equals the absolute value of the determinant of its Kasteleyn matrix,
and then show how this implies that the probability of a pattern in a
random tiling of a board is equal to a minor of the inverse of the 
Kasteleyn matrix of the board. This minor has side proportional to
the number of squares in the pattern. Finally, we will show that
the problem of finding an entry in the inverse of a Kasteleyn matrix
can be reduced to an enumerative problem concerning an Aztec diamond
with one black hole and one white hole.

We will henceforth refer, not to boards, but to their dual graphs, 
which can be seen as having a vertex at the center of every square
and an edge perpendicular to every edge between two squares. For
example, the Aztec diamond will mean for us the object in figure
\ref{fig:azg} and not the object in figure \ref{fig:azb}. This convention
will simplify graph-theoretical arguments considerably.

The results in this section are in part a modern formulation of Kasteleyn's
work, and in part a codification of local folklore, as crafted and passed
down by R. Kenyon, J. Propp, D. Wilson and others. 

\subsection{The Kasteleyn-Wilson matrix}

Consider a finite 
subgraph G of the infinite square lattice with as many white as
black vertices, where the infinite square lattice is colored as
a checkerboard. Let $(v_1,\dotsb ,v_n)$ be its white vertices and 
$(w_1,\dotsb ,w_n)$ its black vertices. (Any ordering from $1$ to $n$
can be chosen.) The
Kasteleyn-Wilson matrix \[K((v_1,\dotsb ,v_n),(w_1,\dotsb ,w_n))\] (or
$K(G)$, by abuse of language) is defined
to be $\determinant{a_{i,j}}_1^n$, where $a_{i,j}$ is
\begin{itemize}
\item $0$, if $\{v_i,w_j\}$ is not an edge of G;
\item $1$, if $\{v_i,w_j\}$ is a horizontal edge of G;
\item $(-1)^k$, if $\{v_i,w_j\}$ is a vertical edge going from row $l$ 
to row $l+1$ and there are $k$ vertices in row $l$ to the left of the edge.
\end{itemize}

We will show that the number of perfect matchings of $G$ is equal to
the absolute value of the determinant of $K(G)$. As could be exprected,
the Kasteleyn-Wilson matrix is only one of several Kasteleyn matrices $K(G)$,
that is, determinants whose absolute values are equal
to the perfect matchings of $G$. We will prove the enumerative property
of the Kasteleyn-Wilson matrix because it is true for any subgraph $G$
of the infinite square lattice, and not only for the Aztec diamond.
Fortunately, the proof can be applied to other Kasteleyn matrices with
minimal cases. In fact, in the last subsection of this section, and in 
following sections, we will use the following convention, which is more
convenient for our purposes than Wilson's:
\begin{itemize}
\item $0$, if $\{v_i,w_j\}$ is not an edge of G;
\item $1$, if $\{v_i,w_j\}$ is a horizontal edge of G;
\item $(-1)^k$, if $\{v_i,w_j\}$ is a vertical edge going from row $l$ 
to row $l+1$ and there are $k$ vertical edges from row $l$ to row $l+1$
to the left of the edge.
\end{itemize}
This convention is valid for the Aztec diamond and for any 
other subgraph $G$ of the infinite square lattice such that any
lattice vertex inside any loop in $G$ is also a vertex of $G$. Thus, only
the material in the following two subsections is valid for any subgraph
of the infinite square lattice. The rest of this work is specific to
the Aztec diamond, although the techniques used are applicable
to other boards.

\subsection{Why does $K$ give us the number of perfect matchings?}

We want to show that the number of perfect matchings of a subgraph $G$ of the
infinite square lattice is equal to the absolute value of the
determinant of $K(G)$. 

We can express $\det(K) $ as
\begin{equation}
\sum_{\pi \in P(\{1,2,...,n\})} \sgn(\pi) \prod_{i=1}^n a_{i,\pi(i)},
\end{equation}

where $P(\{1,2,...,n\})$ is the set of all permutations of $\{1,2,\dotsb n\}$.
Define a map $f$ from the set of all perfect matchings of $G$ to 
$P(\{1,2,...,n\})$ as follows. Any perfect matching can
be expressed in the form
\begin{equation}
\{\{v_1,w_{k(1)}\},\{v_2,w_{k(2)}\},\dotsb ,\{v_n,w_{k(n)}\}\}, 
\end{equation}
where $k$ is a map from $\{ 1,2,\dotsb n\} $ to itself.
The map $f$ takes 
$\{\{v_1,w_{k(1)}\},\dotsb ,\{v_n,w_{k(n)}\}\}$ to $k$. It
is clear that $k$ is a permutation; otherwise there would be unpaired
vertices, as well as vertices belonging to more than one pair. Thus, 
$f$ is well defined. Moreover,
\begin{enumerate}
\item $f$ is injective: If two matchings had the same map $k$, they would be
the same matching.
\item every element $k$ of the image of $f$ satisfies 
$(\prod_{i=1}^n a_{i,k(i)})\ne 0$:
if ${v_i,w_{k(i)}}$ is an edge, then $a_{i,k(i)}$ must be non-zero.
\item if a permutation $\pi $ of ${1,2,...,n}$ satisfies 
$\prod_{i=1}^n a_{i,\pi(i)}\ne 0$,
then, for every $1\le i\le n$, $(v_i,w_{\pi(i)})$ is a valid edge. Moreover,
for $i_1\ne i_2$, $\pi(i_1) \ne \pi(i_2)$, and thus $(v_{i_1},w_{\pi(i_1)})$
 and
$(v_{i_2},w_{\pi(i_2)})$ do not have any vertices in common. Therefore
\begin{equation}\{\{v_1,w_{k(1)}\},\{v_2,w_{k(2)}\},\dotsb ,\{v_n,w_{k(n)}\}\}\end{equation}
is a perfect matching.
\end{enumerate}

Hence $f$ is a one-to-one and onto map from the set of all perfect
matchings of $G$ to the set of all permutations $k$ 
of ${1,2,...,n}$ satisfying
$\prod_{i=1}^n a_{i,k(i)} \ne 0$. Therefore there are as many
perfect matchings as there are non-zero terms in
\begin{equation} \det(K) = 
\sum_{\pi \in P(\{1,2,\dotsb ,n\})} \sgn(\pi) \prod_{i=1}^n a_{i,\pi(i)}.\end{equation}
Every non-zero term is equal to either $1$ or $-1$. 
To prove that the absolute value of
$\det (K)$ equals the number of perfect matchings, we have
to show that all non-zero terms have the same sign. 

\begin{figure}
\centering
	\includegraphics[height=1.5in]{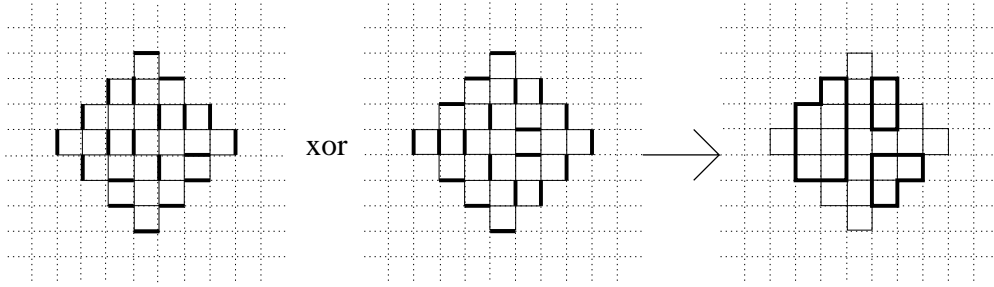}
	\caption{The union of two perfect matchings of the Aztec diamond,
	 minus their intersection}\label{super}
\end{figure}

Let $M = \{\{v_i,w_{\pi(i)}\}\}_{i=1}^n$ and 
$M\prime = \{\{v_i,w_{\pi \prime (i)}\}\}_{i=1}^n$
be two perfect matchings of $G$. By the definition of perfect matching,
every vertex of $G$ is in one edge of $M$ and in one edge of $M\prime $.
It follows that every vertex of $G$ is either in two edges or in no edges
of $(M\cup M\prime)-(M\cap M\prime)$. Therefore 
 $(M\cup M\prime)-(M\cap M\prime)$ consists entirely of loops, that is,
it is a collection of disjoint sets of the form
\begin{equation}\{\{ v_{i_1},w_{j_1}\} ,\{w_{j_1},v_{i_2}\},\{ v_{i_2},w_{j_2}\}, \dotsb
    \{ v_{i_m},w_{j_m}\}, \{w_{j_m},v_{i_1}\}\}.\end{equation}
(See figure \ref{super}.) If a vertex is in two edges of
$(M\cup M\prime)-(M\cap M\prime)$, one of these two edges must be in
$M$, and the other one in $M\prime $. We can assume without loss of generality
that $\{v_{i_1},w_{j_1}\}$ is in $M$, and hence
$\{w_{j_1},v_{i_2}\}$ is in $M\prime $, 
$\{v_{i_2},w_{j_2}\}$ is in $M$, and so on. Then, on one hand,
$j_l = \pi(i_l)$ for $1\leq l\leq m$, and, on the other hand,
$j_l = \pi\prime (i_{l+1})$ for $1\leq l\leq m-1$,
$j_m = \pi\prime (i_1)$. Hence $i_2 = ((\pi\prime)^{-1} \circ \pi)(i_1)$,
$i_3 = ((\pi\prime)^{-1} \circ \pi)(i_2)$,\dots 
$i_1 = ((\pi\prime)^{-1} \circ \pi)(i_m))$.
Thus every loop in $(M\cup M\prime)-(M\cap M\prime)$ induces a cycle in
$(\pi\prime)^{-1} \circ \pi$. It is easy to see that, conversely,
 for every cycle in
$(\pi\prime)^{-1} \circ \pi$ there is a loop in 
$(M\cup M\prime)-(M\cap M\prime)$ that induces it. Because the sign
of a permutation is equal to the product over all its cycles of
$(-1)$ to the power of the length of the cycle minus one, and because
a cycle has length equal to half the number of edges of the loop inducing
it, we have
\begin{equation}\sgn((\pi\prime)^{-1} \circ \pi) = \prod_{\ell\in L} (-1)^{\len(\ell)/2-1},\end{equation}

where $L$ is the set of all loops of $(M\cup M\prime)-(M\cap M\prime)$ and
$\len(\ell ) $ is the number of edges in loop $\ell $. From this equation, from
\begin{equation}\sgn((\pi\prime)^{-1} \circ \pi) = \frac{\sgn(\pi )}{\sgn(\pi \prime )},\end{equation}
and from the fact that all $a_{i,j}$ are $1$ or $(-1)$, it follows that
the result we want to prove in this section, namely,
\begin{equation}\sgn(\pi ) \prod_{i=1}^n a_{i,\pi(i)} = \sgn(\pi \prime)
        \prod_{i=1}^n a_{i,\pi \prime (i)},\end{equation}
is equivalent to
\begin{equation}\prod_{\ell\in L} (-1)^{\len(\ell)/2-1} = \prod_{i=1}^n a_{i,\pi(i)} \cdot
				       \prod_{i=1}^n a_{i,\pi \prime (i)}.\end{equation}
Now, 
\begin{equation}\prod_{i=1}^n a_{i,\pi(i)} \cdot \prod_{i=1}^n a_{i,\pi \prime (i)}
= \prod_{\ell \in L} (\prod_{i\in I(\ell )} a_{i,\pi(i) } \cdot
			\prod_{i\in I(\ell )} a_{i,\pi \prime(i) })\end{equation},
where $I(\ell )$ is the set of indices $i$ of all white vertices in loop
$\ell $. Therefore it is enough for us to prove that
\begin{equation}
(-1)^{\len(\ell)/2-1} =
\prod_{i\in I(\ell )} a_{i,\pi(i) } \cdot
			\prod_{i\in I(\ell )} a_{i,\pi \prime(i) }
\end{equation}
for every loop $\ell $ in $(M\cup M\prime)-(M\cap M\prime)$.

$\prod_{i\in I(\ell )} a_{i,\pi(i) } \cdot
			\prod_{i\in I(\ell )} a_{i,\pi \prime(i) }$
is the product $a_{i,j}$ over all $1\leq i,j\leq n$ such that
$\{v_i,w_j\}$ is an edge of loop $\ell $. Since $a_{i,j}=1$ for
$\{v_i, w_j\}$ horizontal, we can restrict the product to 
$1\leq i,j\leq n$ such that
$\{v_i,w_j\}$ is a vertical edge of loop $\ell $.
What is, specifically, the product of $a_{i,j}$ over all $1\leq i,j\leq n$
such that $\{v_i,w_j\}$ is a vertical edge of loop $\ell $ having
a vertex on row $y_0$ and another in row $y_0+1$, 
where $y_0$ is a given? Let
$\{(x_1,y_0),(x_1,y_0+1)\},\{(x_2,y_0),(x_2,y_0+1)\},\dotsb
\{(x_m,y_0),(x_m,y_0+1)\}$ be all such edges. (We are referring vertices
by their Cartesian coordinates.) Since the horizontal line
$\{(t,y+\frac{1}{2}): t\in (-\infty, \infty)\}$ crosses the loop
an even number of times, $m$ must be even. Let $m = 2\cdot m_0$.
Since all $a_{i,j}$ are $1$ or $(-1)$, we have
\begin{equation}\prod_{j=1}^m a_{J_0(x_j,y_0),J_1(x_j,y_0+1)} =
  \prod_{j=1}^{m_0} \frac{a_{J_0(x_{2j},y_0),J_1(x_{2j},y_0+1)}}
			 {a_{J_0(x_{2j-1},y_0),J_1(x_{2j-1},y_0+1)}} ,\end{equation}
where $J_0(z,w)$ is the index of the white vertex of coordinates
$(z,w)$, and $J_1(z,w)$ is the index of the black vertex of coordinates
$(z,w)$. Since $a_{J_0(x_i,y_0),J_1(x_i,y_0+1)}$ is equal to the
number of vertices of $G$ on row $y_0$ and to the left of $x_{i}$, and
$a_{J_0(x_{i+1},y_0),J_1(x_{i+1},y_0+1)}$ is equal to the
number of vertices of $G$ on row $y_0$ and to the left of $x_{i+1}$,
\begin{equation}\frac{a_{J_0(x_{i+1},y_0),J_1(x_{i+1},y_0+1)}}
			 {a_{J_0(x_{i},y_0),J_1(x_{i},y_0+1)}}\end{equation}
is equal to
the number of vertices of $G$ on row $y_0$ and to the left of $x_{i+1}$
but not of $x_i$. This is the same as the number of vertices of
the infinite square grid on row $y_0$ and to the left of $x_{i+1}$ but
not of $x_{i}$, minus the number of vertices in the grid but not in $G$,
on row $y_0$ and to the left of $x_{i+1}$ but
not of $x_{i}$. The number
of vertices of the grid to the left of $x_{i+1}$ but not of $x_i$ is
equal to the number of squares of the grid 
contained between the edges $((x_i,y),(x_i,y+1))$ and 
$((x_{i+1},y),(x_{i+1},y+1))$. 

It is clear that, when $i$ is even, the squares between 
the edges $((x_i,y),(x_i,y+1))$ and 
$((x_{i+1},y),(x_{i+1},y+1))$ are in the interior of the loop, as are
the vertices on row $y$ of the grid which do not belong to $G$ and which
are the left of $x_{i+1}$ but not
of $x_i$. (Since any vertex on the loop, and, specifically, $x_i$ and
$x_{i+1}$, must be in $G$, we can hencefort refer to these vertices
as ``the vertices on row $y$ of the grid which do not belong to $G$ and
which lie between $x_i$ and $x_{i+1}$''.) Conversely, every square between
rows $y$ and $y+1$ and in the interior of the loop lies between edges
 $((x_i,y),(x_i,y+1))$ and  $((x_{i+1},y),(x_{i+1},y+1))$ for some odd $i$,
and, moreover,
 every vertex on row $y$ of the grid, not in $G$, and in the interior
of the loop lies between $x_i$ and $x_{i+1}$ for some odd $i$.
Hence the number of squares between edges $((x_i,y),(x_i,y+1))$
and  $((x_{i+1},y),(x_{i+1},y+1))$ for $i$ odd equals the number of
squares which are both 
between rows $y$ and $y+1$ and in the interior of the loop, and, furthermore,
the number of vertices not in
$G$ lying on row $y$ between $x_i$ and $x_{i+1}$ for $i$ odd equals the
number of vertices not in $G$ lying on row $y$ and in the interior of the loop.
Thus 
\begin{equation}\prod_{j=1}^m a_{J_0(x_j,y_0),J_1(x_j,y_0+1)} =
  \prod_{j=1}^{m_0} \frac{a_{J_0(x_{2j},y_0),J_1(x_{2j},y_0+1)}}
			 {a_{J_0(x_{2j-1},y_0),J_1(x_{2j-1},y_0+1)}} \end{equation}
equals $(-1)$ to the power of
the number of squares which are both
in the interior of the loop and between rows $y$ and $y+1$
minus the number of such vertices in the interior on row $y$.
Hence, the product of this expression over all rows, that is,
\begin{equation}
\frac{\prod_{i\in I} a_{i,\pi(i)}}{\prod_{i\in I} a_{i,\pi \prime (i)}}
\end{equation}
 equals $(-1)$ to the power of
the number of squares in the interior of the loop minus the number of vertices
lying on row $y$ and in the interior, but not in $G$. 

Therefore 
\begin{equation}
(-1)^{\frac{\len(\ell)}{2}-1} \frac{\prod_{i\in I(\ell )} a_{i,\pi(i)}}
     {\prod_{i\in I(\ell )} a_{i,\pi \prime (i)}}
\end{equation}
is equal to $(-1)$ to the power of the length of the loop, divided by 2,
minus $1$, plus the number of squares inside the loop, minus the number
of vertices in the grid and inside the loop, but not in $G$.
Pick's theorem states that, for polygons whose vertices belong to a square
grid, $A=I+B/2-1,$ where $A$ is the area enclosed by the polygon,
$I$ is the number of grid points inside the polygon, and $B$ is the
number of grid points on the boundary of the polygon. Hence 
our expression is equal to $(-1)$ to the power of the number of
vertices inside the loop, minus the number of vertices in the grid and
inside the loop, but not in $G$. This is the same as the number
of vertices of $G$ inside the loop. Now, the vertices inside the loop
are matched only among themselves, not with the vertices outside the loop,
in either $M$ or $M\prime $. For a set of vertices to be matched only among
themselves, there must be an even number of vertices in the set.
Therefore
the number of vertices of $G$ inside the loop must be even, and
$(-1)$ raised to the power of this number must be $1$. 

Hence
\begin{equation}
(-1)^{\frac{\len(\ell)}{2}-1} \frac{\prod_{i\in I(\ell )} a_{i,\pi(i)}}
     {\prod_{i\in I(\ell )} a_{i,\pi \prime (i)}} = 1
\end{equation}
as we desired to prove. We conclude, by taking the product over all loops
$\ell \in L$, that

\begin{equation}\frac{\sgn(\pi) \prod_{i=1}^n a_{i,\pi(i)}}
       {\sgn(\pi \prime) \prod_{i=1}^n a_{i,\pi \prime(i)}} = 1.\end{equation}

Therefore the terms of
\begin{equation}
\det(K) = \sum_{\pi \in P(\{1,2,...,n\})} \sgn(\pi) \prod_{i=1}^n a_{i,\pi(i)}
\end{equation}
all have the same sign. It follows that the number of perfect matchings, 
that is,
of permutations $\pi$ for which $\prod_{i=1}^n a_{i,\pi(i)}$ is non-zero,
is equal to the absolute value of $\det(K)$.

\subsection{Why does $K^{-1}$ give us the probabilities of patterns?}

Let us have a graph $G$ and a subgraph $H$ whose every vertex has degree $1$.
There is a one-to-one correspondence 
between perfect matchings of $G$ having $H$ as a subgraph and perfect matchings
of $G-H$. (By $G-H$ we mean the graph $(V,E)$, where $V$ is the
set of all vertices of $G$ not in $H$ and $E$ is
the set of all edges of $G$ between two vertices of $G$ not in $H$.)
Therefore the probability that a random perfect matching of $G$ have the
set of edges of $H$ as a subset is equal to the number of perfect matchings
 of $G-H$ divided by the number of perfect matchings of $G$. 

What is, then, the number of perfect matchings of $G-H$? One way to compute it
is to construct a Kasteleyn-Wilson matrix for $G-H$. Another way,
which will soon prove its virtues, is to take the minor of the
Kasteleyn-Wilson matrix for $G$ resulting from the deletion of the rows
and columns corresponding to $H$. Certainly, this minor is not the same
as the Kasteleyn-Wilson matrix for $G-H$; the signs of the matrix entries are
different. Nevertheless, the absolute value
of the determinant of the minor 
is equal to the number of perfect matchings of $G-H$. To prove
this, we need to do the same as in the previous section, namely, show that
all non-zero terms of the expression of the determinant as a
sum over permutations have the same sign.
If we proceed in the same way as before,
we  arrive at a point where the only difference is that
we have $(-1)$ to the power of the number of vertices of $G$ inside a loop
instead of the number of vertices of $G-H$ inside the same loop.
This difference is no difference if there is an even number of vertices
of $H$ inside the loop. Since the vertices of $H$ 
inside the loop cannot be connected with the vertices of $H$ outside 
the loop, the vertices of $H$ inside the loop are paired among themselves.
Thus we have, as we wanted, that
there is an even number of vertices of $H$ inside the loop.  
This is enough for us to show that 
the terms in the expression of the minor of $K(G)$ as a sum over
permutations do not cancel.

Therefore the number of perfect matchings of $G-H$ is equal to the
absolute value of the determinant of the minor of $K(G)$ lacking
the rows and columns corresponding to the vertices of $H$. If we
assign the same probability to every perfect matching of $G$, the probability
that a random perfect matching will have $H$ as a subgraph will be equal to
the number of perfect matchings of $G-H$ divided by the number of
perfect matchings of $G$, that is, to the absolute value of the
determinant of the aforementioned minor of $K(G)$ divided by
the absolute value of the determinant of $K(G)$. By Jacobi's rule,
a corollary of Cramer's rule, this is equal to the absolute value of the
determinant of the minor of $(K(G)^{-1})^T$ consisting of those rows
and columns omitted from the minor in the numerator, that is, of
rows
$1\leq b_1<b_2<\dotsb b_m\leq n$ and columns
$1\leq c_1<c_2<\dotsb c_m\leq n$, where $v_{b_1},v_{b_2},\dotsb v_{b_m}$
and $w_{c_1},w_{c_2},\dotsb w_{c_m}$ are the vertices of $H$.

This is a clear improvement over the
expression $\frac{\|K(G-H)\|}{\|K(G)\|}$. Instead of dealing with determinants
of the size of $G$, we deal with a determinant of the size of $H$.
As we explained in the introduction, we are interested in finding
the probability of small subgraphs, or ``local patterns'', in a large graph 
$G$. We want to find what happens when we have an infinite
sequence of $G$'s whose number of vertices goes to infinity.
Now it is enough for us to examine
a fixed number of entries in each $(K(G_i)^{-1})^T$, and determine their
asymptotic behavior as $i\to \infty.$ 

\subsection{How can we find the entries of $K^{-1}$ for an Aztec diamond
by counting perfect matchings?}

Suppose that we have to compute the determinant of the minor of $K(G)^{-1}$
consisting of rows $1\leq b_1 < b_2 <\dots <b_m\leq n$ and columns
$1\leq c_1 < c_2 <\dots <c_m\leq n$. To do so, we have to compute the
entry $((K(G)^{-1})^T)_{b_i,c_j}$ for $1\leq i,j\leq m$. In other words,
in order to compute the probability of any local pattern, it suffices
to be able to compute an arbitrary entry of the inverse Kasteleyn matrix.
If we have a sequence of graphs $\{G_k\}_{k=1}^\infty$ 
(such as, for example, Aztec diamonds of higher and
higher order) we will know the asymptotic behavior of the probabilities
of local patterns if we know the asymptotic behavior of the entries
in the sequence of matrices $\{K(G_k)^{-1}\}_{k=1}^\infty$.

We can, of course, compute a first minor\footnote{That is, a minor
that has all columns of the matrix of which it is a minor, but one,
and all rows but one.}
 of $K(G)$, and then apply Cramer's
rule, whenever we want an entry of $(K(G)^{-1})^T$.
Unfortunately, it seems very hard
to obtain asymptotic expressions directly from the minors.
We will reduce the problem of computing the minor of $K(G)$
resulting from the deletion of row $i$ and column $j$ to an enumerative
problem whose solution we will be able to represent in a form other than a
determinant.
It would seem, at first sight, that we can use the same
kind of argument we used to answer the previous question, and show
that such a minor is equal (up to sign) to the number of all perfect matchings
of $G$ with white vertex $i$ and black vertex $j$ deleted.
Unfortunately, this is not the case. A first minor without row $i$
and column $j$ is equal to a sum whose number of terms is equal to
number of perfect matchings
of $G$ with white vertex $i$ and black vertex $j$ deleted.       
The problem is that, while every term has absolute value $1$,
not every term has the same sign. It is easy to show, by the same
kind of reasoning we employed in our answer to the first question, that
the matter of whether or not two terms have the same sign can be determined by
examining the loops in the superimposition of the two matchings
corresponding to the two terms. The terms have the same sign
if and only if there is an even number of loops having one of the two deleted
vertices, but not the other, in their interiors. 

\begin{Def}
The {\em Aztec diamond} of order $n$ is a planar graph
consisting of vertices 
\begin{equation}
\{(2r+1,2s) : 0\leq r < n, 0\leq s \leq n\} \cup
\{(2r,2s+1) : 0\leq r \leq n, 0\leq s < n \}
\end{equation}
and of edges
\begin{equation}
\begin{aligned}
\{&\{\{(2r,2s+1),(2r+1,2s+2)\},
   \{(2r+1,2s+2),(2r+2,2s+1)\}, \\
   &\{(2r+2,2s+1),(2r+1,2s)\},
   \{(2r+1,2s),(2r,2s+1)\}\} : 0\leq r< n, 0\leq s < n \}
\end{aligned}
\end{equation}
in Cartesian coordinates.
\end{Def}

\begin{figure}
        \begin{minipage}[b]{0.5\linewidth}
                \centering \includegraphics[height=2in]{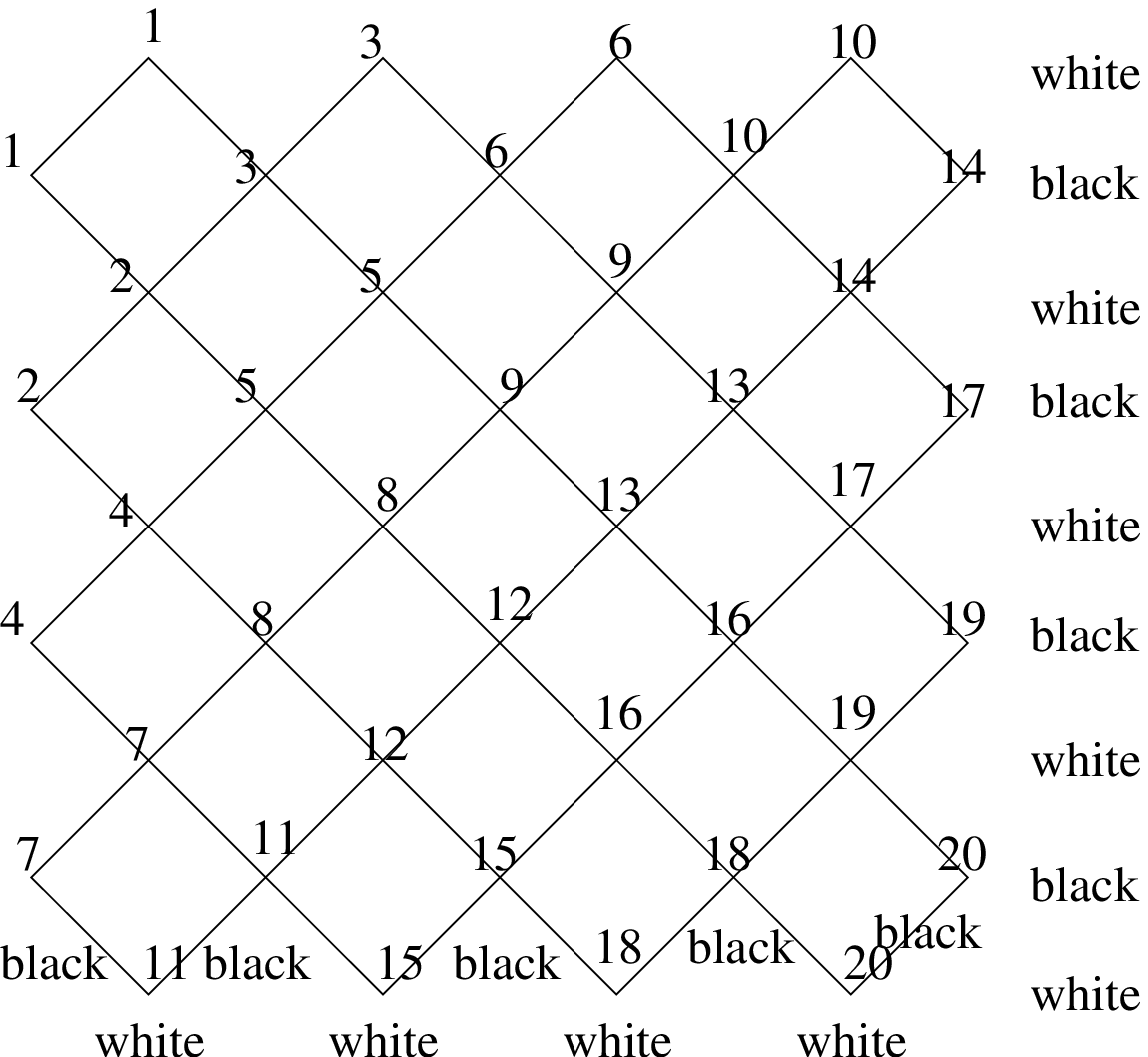}
                \caption{Ordering}\label{fig:indoeu}
        \end{minipage}%
        \begin{minipage}[b]{0.5\linewidth}
                \includegraphics[height=2in]{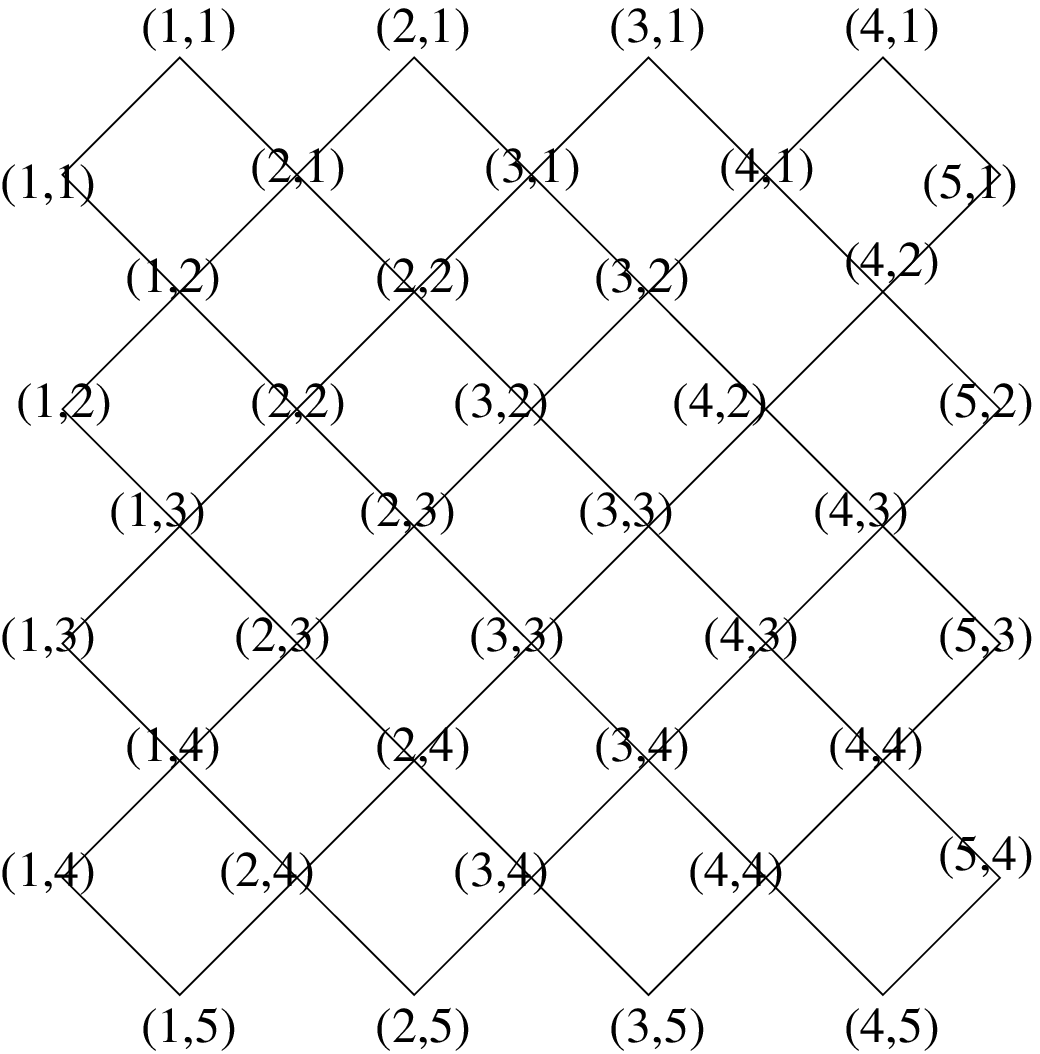}
                \caption{Coordinates} \label{fig:coor}   
        \end{minipage}%
\end{figure}

It will soon become apparent that, for, our purposes, the
system of coordinates in Figure \ref{fig:coor} is more convenient than
Cartesian coordinates. It will also become clear why we draw
the Aztec diamond as if on an infinite square grid tilted $45$ degrees
from the ``natural'' direction. For now, let us notice that, if we
color the Aztec diamond as a checkerboard, all vertices on a column
have the same color, as do all vertices on a row. Let us also
use the system of coordinates in Figure 
\ref{fig:coor} instead of the Cartesian 
system, and refer to the edge consisting of white vertex $(x_0,y_0)$
and black vertex $(x_1,y_1)$ as $((x_0,y_0),(x_1,y_1))$.

We can reformulate the rule for comparing terms' signs so
that it does not mention loops.

\begin{figure}
\centering
	\includegraphics{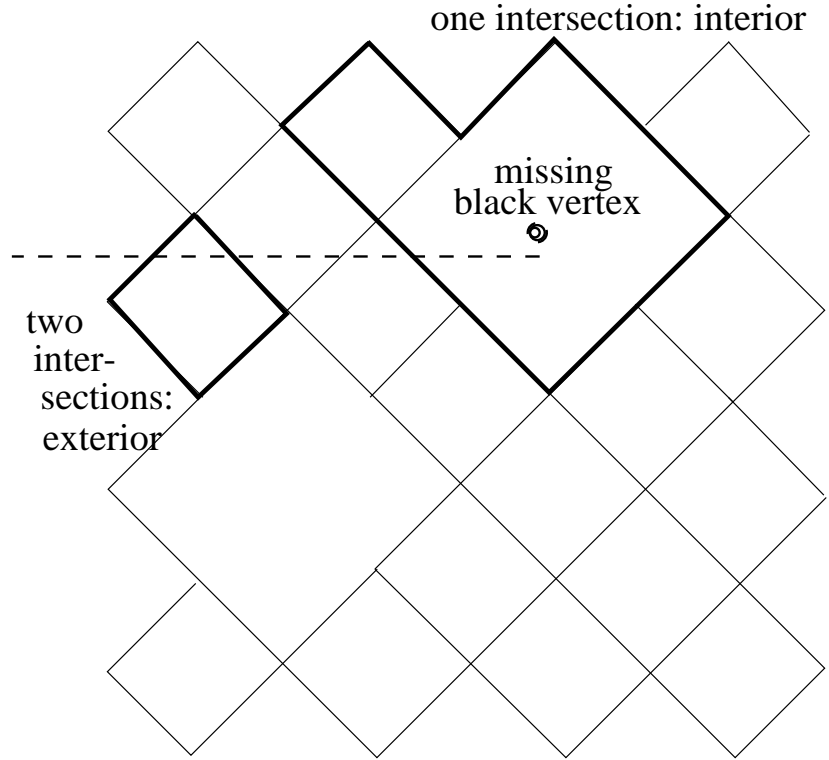}
	\caption{Counting intersections}\label{fig:whynot}
\end{figure}

\begin{Lem}\label{Lem:opposum}
Let $A$ and $B$ be two perfect matchings of the Aztec diamond of 
order $n$ with the white vertex at $(w_0,w_1+d_1)$ and the black
vertex at $(w_0+d_0,w_1)$ deleted\footnote{We shall henceforth use
the system of coordinates in Figure
\ref{fig:coor} instead of the Cartesian system.}.
The following two conditions are equivalent:
\begin{enumerate}
\item $A\cup B - A\cap B$ has an even number of loops
containing exactly one of the two deleted vertices.
\item $w(A) \cong w(B) \mod 2$, where 
$w(T)$ is the number of edges of the form $((i-1,w_1+1),(i,w_1))$,
$1<i<w_0+d_0$; $((i,w_1+1),(i,w_1))$, $1\geq i<w_0+d_0$;
$((i,w_1+d_1),(i,w_1+d_1-1)$, $1\geq i<w_0$, and
$((i,w_1+d_1),(i+1,w_1+d_1-1))$, $1\geq i<w_0$ in a perfect matching $T$.
\end{enumerate}
\end{Lem}
\begin{proof}
The sum $w(A) + w(B)$ is congruent, modulo $2$, to
the total number of edges in $A\cup B - A\cap B$ consisting 
of a black vertex $(x,y)\in D$
and white vertex $(x-1,y+1)$ or $(x,y+1)$. Given a loop $\ell $
in $A\cup B - A\cap B$, the number of edges in it consisting
of a black vertex $(x,y)\in E$
and white vertex $(x-1,y+1)$ or $(x,y+1)$ is equal to the number
of times the loop crosses a ray with its end slightly below
the deleted black vertex and with diagonal direction with respect
to the square grid. (See figure \ref{fig:whynot}.)
 This number is even if  the deleted
black vertex is in the exterior of the loop, and odd if it is in the interior.
The same holds for $F$ and the deleted white vertex. Hence
a loop has an even number of edges consisting of a black vertex $(x,y)\in D$
and white vertex $(x-1,y+1)$ or $(x,y+1)$ if and only if it has exactly
one of the two deleted vertices in its interior. We conclude the proof
by summing over all loops.
\end{proof}

It follows that $|\sum (-1)^{w(T)}| = 
|\sum_{\pi } \sgn(\pi ) \prod_{i=1}^n a_{i,\pi(i)} |$, where
the sum on the left is over all perfect matchings $T$ of the Aztec
diamond with white vertex at $(w_0,w_1+d_1)$ and black
vertex at $(w_0+d_0,w_1)$, and $\determinant{a_{i,j}}_1^n$ is the
first minor of the Kasteleyn matrix $K(G)$ of the (intact) Aztec diamond $G$
where the row corresponding to the white vertex at $(w_0,w_1+d_1)$
and the column corresponding to the black vertex at $(w_0+d_0,w_1)$.
Thus, if we compute $|\sum (-1)^{w(T)}|$ (and this is essentially
an enumerative problem, which we shall solve enumeratively), we will
know the absolute value of 
$\sum_{\pi } \sgn(\pi ) \prod_{i=1}^n a_{i,\pi(i)}$, which
will give us the absolute value of the entry of $(K(G)^{-1})^T$ to
be computed, but not its sign. We will find the sign now.

It is better, for this particular goal, to
express the entry of the inverse Kasteleyn matrix of the Aztec diamond, 
not as the ratio of the determinant of the minor of the Kasteleyn matrix
resulting from the deletion of 
column $i_0$ and row $j_0$, divided by the determinant of the Kasteleyn
matrix, multiplied by $(-1)^{i_0+j_0}$, but rather as the ratio
of the determinant $\determinant{b_{i,j}}_1^{n-1}$ to the determinant of the
Kasteleyn matrix $K(G) = \determinant{a_{i,j}}_1^n$, where $b_{i_0,j}=0$
for $j\ne j_0$, $b_{i,j_0}=0$ for $i\ne i_0$, $b_{i_0,j_0}=1$,
$b_{i,j}=a_{i,j}$ for $i\ne i_0$, $j\ne j_0$. We can then ask whether,
for a permutation $\pi $ of $\{1,2,\dotsb ,n-1\}$, 
$\sgn(\pi ) \cdot \prod_{i} b_{i,\pi(i)}$
has the same sign as $(-1)^{w(T)}$, where $T$ is the perfect matching 
corresponding to $\pi$. (The answer will be the same for all permutations, 
so we have to ask it for only one permutation (or matching) we choose.)
If the signs are the same, then
$\sum (-1)^{w(T)} = 
(((K(G))^{-1})^T)_{i_0,j_0} \determinant{K(G)}_1^n$; if
the sign are different, then
$\sum (-1)^{w(T)} = 
- (((K(G))^{-1})^T)_{i_0,j_0} \determinant{K(G)}_1^n$.

For our search for the sign, we need to fix an ordering of the black and
white vertices of the Aztec diamond. The ordering in Figure
\ref{fig:indoeu} will prove itself convenient.

We want to prove that, for any $(w_0,w_1,d_0,d_1)$, 
\[
\sum (-1)^{w(T)} = (-1)^{d_0+d_1+1}\cdot 
\sum_{\pi } \sgn(\pi ) \prod_{i=1}^n a_{i,\pi(i)}. \label{eq:mandrake}
\]

The most straightforward method of proof consists of comparing the signs
of $\sum (-1)^{w(T)}$ and $\sgn(\pi ) \prod_{i=1}^n a_{i,\pi(i)}$ for
some tiling $T$ and its corresponding permutation $\pi $. Unfortunately,
this is also quite cumbersome, as which tilings are possible depends
on the order of $w_0$, $w_0+d_0$, $w_1$ and $w_1+d_1$. In order not to
bore the reader with twenty-four different cases, we will proceed
explicitly only for the six cases corresponding to $d_0,d_1>0$, which
can be treated as four different cases.

\subsubsection{Case 1: $w_1\leq w_0<w_1+d_1$}

We choose to examine the matching having $/$-edges at
$((i,w_1),(i,w_1))$, $w_1\leq i< w_0+d_0$, $((w_1,j),(w_1,j))$,
$w_1< j < w_1+d_1$, and $((i,w_1+d_1),(i+1,w_1+d_1-1))$, 
$w_1\leq i < w_0$, where $((x_0,y_0),(x_1,y_1))$ denotes the
edge consiting of a white vertex at $(x_0,y_0)$ and a black vertex
at $(x_1,y_1)$. All other vertices are covered by $\setminus$-edges:
$((x,y),(x+1,y))$ for $x\geq y$, given that at most one of the two
conditions $w_1\leq x <w_0+d_0$, $y = w_1$ is fulfilled, and
$((x,y),(x,y-1))$ for $x<y$, given that at most one of 
$x = w_1$, $w_1< y <w_1+d_1$ holds, as well as at most one of
$w_1\leq x < w_0$, $y = w_1+d_1$. (Note that we denote by
$/$-edges and $\setminus$-edges what we called ``vertical edges'' and
``horizontal edges'' before. This is just a change in notation with 
the purpose that the fact that we now draw the Aztec diamond tilted
$45$ degrees will not confuse the reader.)
It is easy to verify that this
is a perfect matching. For this matching, $(-1)^{w(T)}$ is 
$(-1)^{w_0+w_1}$. 

Each matching corresponds to a permutation of $\{1,2,\dotsb , n\}$.
(We have permutations of $\{1,2,\dotsb , n\}$, and not of 
$\{1,2,\dotsb , n-1\}$, because we are working with 
$\determinant{b_{i,j}}_1^n$, and not with a first minor of
$K(G)=\determinant{a_{i,j}}_1^n$. The convenience of this choice,
which we made without justification, will be clear by the end of this
paragraph.) The permutation $\pi$ corresponding
to the matching we are now examining consists of one cycle.
The length of this cycle is equal to the number of vertical lozenges, plus 
one.
(This is not true in general. However, it is
true in cases similar to the one we are currently examining, in which
white vertices $\{ f_1,f_2,\dotsb f_m\} $ and
black vertices $\{ g_1,g_2,\dotsb g_m\} $, which would be paired
in the form $(f_i,g_i)$, $1\leq i\leq m$ in the all-$\setminus$  matching, are
paired in the form $(f_i,g_{i+1})$, $1\leq i<m$, and vertices
$f_m$ and $g_1$ are missing.) Hence $\sgn(\pi)$ is equal to $(-1)$
to the power of the number of vertical lozenges, that is, 
$2(w_0-w_1)+d_0+d_1-1$ lozenges.

We now need only examine $\prod_{i=1}^m b_{i,\pi(i)}$. Among the edges
in our perfect matching, only $((i,w_1+d_1),(i+1,w_1+d_1-1))$,
$w_1\leq i<w_0$, correspond to entries equal to $(-1)$ in the determinant
$\determinant{b_{i,j}}$. Therefore $\prod_{i=1}^m b_{i,\pi(i)} = (-1)^{w_0-w_1}$.

We conclude that our matching $T$, for which 
$(-1)^{w(T)} = (-1)^{w_0+w_1-1}$, contributes 
$(-1)^{2(w_0-w_1)+d_0+d_1-1}\cdot (-1)^{w_0-w_1}$, to the determinant.
Therefore an arbitrary matching $T$ contributes
\begin{equation}
(-1)^{w(T)-(w_0+w_1-1)}\cdot 
(-1)^{2(w_0-w_1)+d_0+d_1-1}\cdot (-1)^{w_0-w_1},
\end{equation}
that is,
\begin{equation}
(-1)^{w(T)} (-1)^{d_0+d_1+1},
\end{equation}
 to the determinant. Therefore
\[
\sum (-1)^{w(T)} = (-1)^{d_0+d_1+1}\cdot 
\sum_{\pi } \sgn(\pi ) \prod_{i=1}^n a_{i,\pi(i)} \] for case 1.

\subsubsection{Case 2: $w_0\geq w_1+d_1$}

We choose to examine the matching having $//$-edges at 
$((i,w_1),(i,w_1))$, $w_1<i<w_0+d_0$, at
$((w_0+1,j),(w_0,j+1))$, $w_1+d_1\leq j\leq w_0$, at
$((w_1,j),(w_1,j))$, $w_1\leq j\leq w_0$, and at
$((i+1,w_0),(i,w_0+1))$, $w_1\leq i<w_0$. All other vertices are covered
by $\setminus$-edges. The $(-1)^{w(T)}$ of this matching is 
$(-1)^{d_1+1}$.
The sign of the permutation is 
\begin{equation}
(-1)^{(w_0+d_0-w_1-1)+(w_0-(w_1+d_1)+1)+(w_0-w_1+1)+(w_0-w_1)},
\end{equation}
that is, $(-1)^{d_0+d_1+1}$. The product $\prod_{i=1}^m b_{i,\pi(i)}$
is $(-1)^{(w_0-(w_1+d_1)+1) + (w_0-w_1)}$, that is, $(-1)^{d_0+d_1}$.
Hence each matching $T$ 
contributes 
\begin{equation}
(-1)^{w(T)} (-1)^{d_0+d_1+1}
\end{equation} to the determinant.
Therefore
\[
\sum (-1)^{w(T)} = (-1)^{d_0+d_1+1}\cdot 
\sum_{\pi } \sgn(\pi ) \prod_{i=1}^n a_{i,\pi(i)} \] for case 2.
\subsubsection{Case 3: $w_0<w_1<w_0+d_0$}

We choose to examine the matching having $//$-edges at
$((w_0,j),(w_0,j))$, $w_0\leq j<w_1+d_1$, at
$((i,w_0),(i,w_0))$, $w_0<i\leq w_0+d_0$, and at
$((w_0+d_0,j+1),(w_0+d_0+1,j))$, $w_0\leq j<w_1$. All other vertices are 
covered by $\setminus$-edges. The $(-1)^{w(T)}$ of this matching is
$(-1)^{w_0+w_1+1}$. The sign of the permutation is
$(-1)^{d_0+d_1}$. The product $\prod_{i=1}^m b_{i,\pi(i)}$
is $(-1)^{w_1-w_0}$. Therefore each matching $T$ 
contributes 
\begin{equation}
(-1)^{w(T)}\cdot (-1)^{d_0+d_1+1}
\end{equation}
to the determinant.
Therefore
\[
\sum (-1)^{w(T)} = (-1)^{d_0+d_1+1}\cdot 
\sum_{\pi } \sgn(\pi ) \prod_{i=1}^n a_{i,\pi(i)} \] for case 3.
\subsubsection{Case 4: $w_1\geq w_0+d_0$}

We choose to examine the matching having $//$-edges at
$((w_0,j),(w_0,j))$, $w_0\leq j<w_1+d_1$, at
$((i,w_0),(i,w_0))$, $w_0<i\leq w_1$, at
$((w_1,j+1),(w_1+1,j))$, $w_0\leq j\leq w_1$, and at
$((i,w_1+1),(i+1,w_1))$, $w_0+d_0\leq i<w_1$.
All other vertices are covered by $\setminus$-edges. 
The $(-1)^{w(T)}$ of this matching is
$(-1)^{d_0+1}$. The sign of the
permutation is $(-1)^{d_0+d_1+1}$.
The product 
$\prod_{i=1}^m b_{i,\pi(i)}$  is $(-1)^{d_0+1}$. Therefore each matching
$T$ contributes \begin{equation} (-1)^{w(T)} (-1)^{d_0+d_1+1}\end{equation}
to the determinant.
Therefore
\[ 
\sum (-1)^{w(T)} = (-1)^{d_0+d_1+1}\cdot 
\sum_{\pi } \sgn(\pi ) \prod_{i=1}^n a_{i,\pi(i)} \] for case 4.

\section{An enumerative problem}

Consider an
Aztec diamond of side, or order, $n$, with black vertex $(w_0+d_0,w_1)$
and white vertex $(w_0,w_1+d_1)$ missing. Our task in this section is
 to compute \begin{equation} \sum_T (-1)^{w(T)} \end{equation}
where $T$ ranges over all perfect matchings of this Aztec diamond
with two missing vertices.
Enumeratively speaking, we must count the number of perfect
matchings, counting as ``negative'' any matching $T$
for which $(-1)^{w(T)}=(-1)$. The function $w(T)$, 
as defined in the previous section,
gives, for a perfect matching
$T$, 
the number of edges consisting of a black vertex $(x,y)\in D$
and white vertex $(x-1,y+1)$ or $(x,y+1)$. $D$ is the union of two
set of black vertices: $E = \{(i,w_1) : i<w_0+d_0\}$ and 
$F = \{(i,w_1+d_1) : i\leq w_0\}$.

Throughout this section, we will assume $d_0>0$, $d_1>0$. At the end it
will become clear that this assumption involves no loss of generality.
One of our tools for attacking the special case $d_0>0$, $d_1>0$ will be
the EKLP Lemma, a classical result in \cite{EKLP} which we will
soon explain. First of all, we must define two kinds of subsets of the 
Aztec diamond. They will be our intermediate objects of study.

\begin{Def} An {\em $n\times m$ black-edged Aztec rectangle} with dents at 
$1\leq x_1<x_2<...<x_m\leq n+1$ is a graph consisting of white vertices 
$(i,j),$ $1\leq i\leq n,$ $1\leq j\leq m$, 
black vertices $(i,j),$ $1\leq i\leq n+1$, $1\leq j\leq m-1$, 
and black vertices $(i,m)$ such that $\forall 1\leq k\leq m$ we have 
$i\ne x_k$.
      An {\em $n\times m$ white-edged Aztec rectangle} with teeth at
$1\leq y_1<y_2<...<y_m\leq n$ is a graph consisting of white vertices $(i,j),$
$1\leq i\leq n,$ $1\leq j\leq m,$ black vertices $(i,j),$ 
$1\leq i\leq n+1,$ $1\leq j\leq m,$ and
white vertices $(i,m+1)$ s.t. $\exists 1\leq k\leq m$ with $i=x_k$.
\end{Def}

Let the number of perfect matchings of an $n\times m$ black-edged 
Aztec rectangle with dents
at $1\leq x_1<x_2<...<x_m\leq n+1$ be called 
$D_{n,m}(x_1,x_2,...,x_m).$
What is the number of perfect matchings of an $n\times m$ white-edged Aztec 
rectangle $Y$ with teeth
at $1\leq y_1<y_2<...<y_m\leq n$? Every white vertex (or ``tooth'') 
$(y_i,m+1)$ must
be covered by a horizontal or vertical edge, which will also cover
black vertex $(y,m)$ or $(y+1,m)$, respectively. Thus, any matching of $Y$ is 
composed of 
\begin{itemize}
\item a matching of a black-edged $n\times m$ 
Aztec rectangle with dents at
$1\leq x_1<x_2<...<x_m\leq n+1,$ where 
$x_i-y_i$ is either $0$ or $1$;
\item the unique matching of the remaining region.
\end{itemize}
 Therefore the number of matchings of $Y$ is the $m$-fold sum
\begin{equation}
E_{n,m}(y_1,...,y_m)=\sum_{x_i=y_i \text{ or } y_i+1} D_{n,m}(x_1,...,x_m)
\label{eq:2} \end{equation}

where we assume that $D_{n,m}(x_1,x_2,...,x_m)=0$ when any two $x_i$'s are 
equal. 

In the same way that, given $D_{n,m},$ we have found $E_{n,m},$ given
$E_{n,m},$ we can find $D_{n,m+1}$. We have again an $m$-fold sum:

\begin{equation}
D_{n,m+1}(x_1,...,x_{m+1})=\sum_{x_i\leq y_i<x_{i+1}} E_{n,m}(y_1,...,y_m)
\label{eq:3}
\end{equation}

where we assume that $E_{n,m}(x_1,x_2,...,x_m)=0$ when any two $x_i$'s are 
equal. 

We will now be able to prove the following by induction. The proof
in \cite{EKLP} does not use Lemma \ref{Lem:A} explicitly.

\begin{Lem}[EKLP Lemma]\label{Lem:eklp}
The number of matchings of an $n\times m$ black-edged Aztec rectangle with
dents $1\leq x_1<x_2<\dotsb <x_m\leq n+1$ is

\begin{equation}
D_{n,m}(x_1,x_2,...,x_m) = \frac{2^{\frac{m(m-1)}{2}}}{(m-1)!!} 
			\determinant{x_i^{j-1}}_1^m,
\label{eq:DE} \end{equation}

where $k!! = 1!\, 2!\, 3!\, \dotsb (k-1)!\, k!$. 
The number of matchings of an $n\times m$ white-edged Aztec rectangle with dents
at $1\leq y_1<y_2<\dotsb <y_m\leq n$ is

\begin{equation}
E_{n,m}(y_1,y_2,...,y_m) = \frac{2^{\frac{m(m+1)}{2}}}{(m-1)!!} 
			\determinant{y_i^{j-1}}_1^m
\label{eq:ED} \end{equation}
\end{Lem}

For the proof of this lemma, we need a special case of another lemma that
we will later state in its full generality.

\begin{Lem}[Lemma A (special case)]\label{Lem:A}
Let $A$ be an operator carrying polynomials to
polynomials of the same or lesser degree:
\begin{equation}
\begin{matrix}
x^0 & \to & & & a_{0,0} x^0\\
x^1 & \to & & a_{1,1} x^1 + & a_{1,0} x^0\\
x^2 & \to & a_{2,2} x^2 + & a_{2,1} x^1 + & a_{2,0} x^0\\
\dotsb
\end{matrix}
\end{equation}

Then \begin{equation}
\determinant{A(x^{j-1})(x_i)}_1^m = 
 a_{0,0}\cdot a_{1,1} \dotsb a_{m-1,m-1} \cdot \determinant{x_i^{j-1}}_1^m,
\label{eq:no}
\end{equation}
where 
\begin{equation} 
A(x^j)(x_i) = 
      a_{j,j} x_i^j + a_{j,j-1} x_i^{j-1} + \dotsb + a_{j,0} x_i^0.
\end{equation}
\end{Lem}
\begin{proof}
The determinant on the left side of (\ref{eq:no}) can be obtained from
the determinant on the right side by elementary column operations.
\end{proof}

\begin{proof}[Proof of EKLP Lemma]
The base case $(m=1)$ is trivial. The inductive step  $D_{n,m}\to E_{n,m}$
follows directly from (\ref{eq:3}) and from Lemma A (special case).

The inductive step $E_{n,m}\to D_{n,{m+1}}$ requires some more work.

\begin{equation}
\begin{aligned}
D_{n,m+1}(x_1,\dotsb x_{m+1}) &= 
   \sum_{x_i\leq y_i<x_{i+1}} E_{n,m}(y_1,...,y_m) \\
  &= \sum_{x_i\leq y_i<x_{i+1}} \frac{2^{\frac{m(m+1)}{2}}}{(m-1)!!} 
	\determinant{y_i^{j-1}}_1^m \\
  &= \frac{2^{\frac{m(m+1)}{2}}}{(m-1)!!} 
	\determinant{x_i^{j-1}+(x_i+1)^{j-1}+\dotsb (x_{i+1}-1)^{j-1}}_1^m \\
  &= \frac{2^{\frac{m(m+1)}{2}}}{(m-1)!!} 
	\determinant{\frac{1}{j} (B_j(x_{i+1}) - B_j(x_i))}_1^m,
\end{aligned}
\end{equation}
where $B_r(y)$ is the Bernoulli polynomial of degree $r$, 
which has the property $\sum_{n=M}^{N-1} n^q = (B_{q+1}(N)-B_{q+1}(M))/(q+1)$.
For a definition of $B_r(y)$, consult \cite{Rademacher}. The only
property of $B_r(y)$
we need to know is that it is a polynomial of degree r whose
leading coefficient is one. We will adopt, for the sake of convenience,
the convention $B_0(y) = y^0$. We can now continue:
\begin{align}
D_{n,m+1}(x_1,\dotsb x_{m+1})
  &= \frac{2^{\frac{m(m+1)}{2}}}{m!!}
	\determinant{B_j(x_{i+1}) - B_j(x_i)}_1^m \\
  &= \frac{2^{\frac{m(m+1)}{2}}}{m!!} 
	\begin{vmatrix}
	1 & B_1(x_1) & B_2(x^1) & \dotsb & B_m(x_1) \\
	1 & B_1(x_2) & B_2(x_2) &\dotsb & B_m(x_2)  \\
	\dotsb \\
	1 & B_1(x_{m+1}) & B_2(x_{m+1}) &\dotsb & B_m(x_{m+1}) 
	\end{vmatrix},\\
\intertext{where we have pasted a new column onto the left edge of the
determinant and a new row onto the top edge, and then added the first
row to the second one, the second to the third one, and so on, in succession}
  &= \frac{2^{\frac{m(m+1)}{2}}}{m!!} 
	\determinant{x_i^{j-1}}_1^{m+1},
\end{align}
by Lemma A (special case).
\end{proof}

Let us now count matchings with weights $1$ and $(-1)$: a matching may count
as one matching or as minus one matchings. Consider an $n\times m$ black-edged
Aztec rectangle with dents at $1\leq x_1<x_2<\dotsb x_{m-1}\leq n+1$,
and, in addition, a dent at $w_0$. Let us multiply the total number of 
matchings by $(-1)$ if there is an odd number of $x_i$'s smaller than $w_0$. 
(Here we are counting either all matchings as positive matchings or all
as negative matchings.) What is, then, this total, weighted number of
matchings? Conveniently, it is
\begin{equation} 
 \frac{2^{\frac{m(m-1)}{2}}}{(m-1)!!} 
  \begin{vmatrix}
	w_0^{j-1} \\
	x_1^{j-1} \\
	x_2^{j-1} \\
	\dotsb \\
	x_{m-1}^{j-1}
  \end{vmatrix}
\end{equation}

We have just taken (\ref{eq:DE}) and shifted the $w_0^j$ row corresponding
to the dent at $w_0$ as many positions upwards as there are $x_i$'s smaller
than $w_0$.

What is the weighted number of matchings of a white-edged $n\times m$ Aztec
rectangle $Y$ with teeth at $y_1,y_2,\dotsb y_{m-1}$ and a black hole
at $(w_0,m)$? (The weight of each matching is  $(-1)$ to the power of the number
of its edges consist of a black vertex $(i,w_1)$ and white vertex
$(i-1,w_1+1)$ or $(i,w_1+1)$, where $1\leq i<w_0+d_0$.)
A cursory examination makes clear 
that a matching of $Y$ will have weight $1$ if and only if the black-edged 
$n\times m$
Aztec rectangle the matching outlines has an even number of dents with indices
lower than $w$. Thus the weighted number of matchings is
\begin{equation}
\sum_{x_i=y_i or y_i+1}
\frac{2^{\frac{m(m-1)}{2}}}{(m-1)!!} 
  \begin{vmatrix}
	w_0^{j-1} \\
	x_1^{j-1} \\
	x_2^{j-1} \\
	\dotsb \\
	x_{m-1}^{j-1}
  \end{vmatrix},
\end{equation}
which is the same as
\begin{equation}
\frac{2^{\frac{m(m-1)}{2}}}{(m-1)!!} 
\begin{vmatrix}
	w_0^{j-1} \\
	y_1^{j-1} + (y_1+1)^{j-1} \\
	y_2^{j-1} + (y_2+1)^{j-1} \\
	\dotsb \\
	y_{m-1}^{j-1} + (y_{m-1}+1)^{j-1}
\end{vmatrix}.
\end{equation}

Here we face a difficulty. We would like to use Lemma \ref{Lem:A}. 
Unfortunately, we have $w_0^{j-1}$ on the top row. What we need is to 
find something which
is transformed into $w_0^{j-1}$ by the linear operations taking a row of
the form $x^{j-1}$ to $x^{j-1} + (x+1)^{j-1}$. First of all, we need a 
stronger version of Lemma \ref{Lem:A}.

\begin{Lem}[Lemma A, stronger version]
Let A be an operator carrying polynomials to
polynomials of the same or lesser degree:

\begin{equation}
\begin{matrix}
x^0 & \to & & & a_{0,0} x^0\\
x^1 & \to & & a_{1,1} x^1 + & a_{1,0} x^0\\
x^2 & \to & a_{2,2} x^2 + & a_{2,1} x^1 + & a_{2,0} x^0\\
\dotsb
\end{matrix}
\end{equation}

Then \begin{equation}
\determinant{\sum_{1\leq k\leq l_i} b_{k,i} A(x^{j-1})(x_{k,i})}_1^m 
= a_{0,0} ... a_{m-1,m-1} 
     \determinant{\sum_{1\leq k\leq l_i} b_{k,i} x_{k,i}^{j-1}}_1^m,
\end{equation}
                                                                         
for any $l_i$, $b_{k,i},$ $1\leq k\leq l_i,$ $1\leq i\leq m$, where 
\begin{equation}
A(x^j)(x_i) = 
      a_{j,j} x_i^j + a_{j,j-1} x_i^{j-1} + \dotsb + a_{j,0} x_i^0.
\end{equation}
\end{Lem}
\begin{proof}
The special case followed from the fact that
the same sequence of column operations transforms the row 
$A(x^{j-1})(y)$ into $y^j$
and the row $A(x^j)(z)$ into $z^j$, for any values $y$,$z$. Therefore the same
sequence transforms the row $a\cdot A(x^j)(y)+b\cdot A(x^j)(z)$ 
into $a\cdot y^j+b\cdot z^j$.
\end{proof}

We need $a_k$'s such that 
$\sum_k a_k A(x^j)(v_k) = w_0^j$ for $0\leq j<m$, 
where $A$ is the operator taking $x^j$ to $x^j + (x+1)^j$. 

\begin{Def}
$\Delta $ is the operator taking $x^j$ to $(x+1)^j-x^j$.
\end{Def}

Then $A = (2I+\Delta)$, where $I$ is the identity operator, and we have,
formally,

\begin{equation}
\begin{aligned}
(2I+\Delta)^{-1} &= \frac{1}{2} \cdot (I+\frac{\Delta }{2})^{-1} \\
                 &= \frac{1}{2} \cdot \sum_{j=0}^\infty (-1)^j 
						\frac{\Delta }{2}^j
\end{aligned}
\end{equation}

Since $\Delta^j$ vanishes on polynomials of degree smaller than $j$,
and we are dealing with polynomials of degree at most $m-1$, we can use
just the first $m$ terms of the series:

\begin{equation}
(1/2) \sum_{j=0}^{m-1} (-1)^j (\Delta/2)^j
\end{equation}

The reader can easily check that, for any $x^k$ with $k<m$,
\begin{equation}
(2I+\Delta) ((\frac{1}{2} \sum_{j=0}^{m-1} (-1)^j 
			(\frac{\Delta}{2})^j) (x^k))
\end{equation}
gives $x^k$ plus  a constant times $\Delta^m x^k$, and, since $m>k$,
$\Delta^m x^k$ is zero. One can also check the same for
\begin{equation}
((\frac{1}{2} \sum_{j=0}^{m-1} (-1)^j (\frac{\Delta}{2})^j) 
						(x^k)) (2I+\Delta).
\end{equation}
It is quite convenient that $(2I+\Delta )$ and 
\begin{equation}
((\frac{1}{2}  \sum_{j=0}^{m-1} (-1)^j (\frac{\Delta}{2})^j) (x^k))
\end{equation}
commute, and that their composition is the same as the identity operator
for the domain we are interested in. We will use
the shorthand $(2I+\Delta)^{-1}$ for
\begin{equation}
((\frac{1}{2} \sum_{j=0}^{m-1} (-1)^j (\Delta/2)^j) (x^k))
\end{equation}
without any compuctions.

If we define $a_0,a_1,\dotsb, a_{m-1}$ by
\begin{equation}
(\frac{1}{2} \sum_{j=0}^{m-1} (-1)^j (\Delta/2)^j)(f(x)) = 
                                 \sum_{i=0}^{m-1} a_i f(x+i),
\end{equation}

then 
\begin{equation}
\begin{aligned}
\sum_{i=0}^{m-1} a_i ((2I+\Delta)(x^j)(w_0+i)) &=
     (((\frac{1}{2} 
        \sum_{j=0}^{m-1} (-1)^j (\Delta/2)^j)(2I+\Delta))(x^j))(w_0) \\
&= (x^j)(w_0) \\
&= w_0^j,
\end{aligned}
\end{equation}
for $0\leq j<m,$ as we desired.

Hence, by Lemma A,

\begin{equation}
\begin{vmatrix}
w_0^{j-1} \\
y_1^{j-1} + (y_1+1)^{j-1} \\
y_2^{j-1} + (y_2+1)^{j-1} \\
\dotsb \\
y_{m-1}^{j-1} + (y_{m-1} + 1)^{j-1}
\end{vmatrix}
\end{equation}

is equal to
\begin{equation}
2^m\cdot 
\begin{vmatrix}
(\frac{1}{2} \sum_{k=0}^{m-1} (-1)^k (\frac{\Delta }{2})^k)(x^{j-1})(w_0) \\
y_1^{j-1} \\
y_2^{j-1} \\
\dotsb \\
y_{m-1}^{j-1}
\end{vmatrix},
\end{equation}

or, in shorthand,

\begin{equation}
2^m\cdot
\begin{vmatrix}
((2I+\Delta)^{-1} (x^{j-1}))(w_0) \\
y_1^{j-1} \\
y_2^{j-1} \\
\dotsb \\
y_{m-1}^{j-1}
\end{vmatrix}
\end{equation}

Now, for every polynomial $p$, 
\begin{equation}
p(w_0) = (((I+\Delta)^{w_0-1})(p))(1).
\end{equation}
Hence we can write

\begin{equation}
2^m\cdot
\begin{vmatrix}
((I+\Delta)^{w_0-1} ((2I+\Delta)^{-1} (x^{j-1})))(1) \\
y_1^{j-1} \\
y_2^{j-1} \\
\dotsb \\
y_{m-1}^{j-1}
\end{vmatrix}
\end{equation}

and eliminate all terms of degree $m$ or higher. Hence the weighted number
of matchings of a white-edged $n\times m$ Aztec rectangle with teeth at
$y_1,y_2,\dotsb y_{m-1}$ and a black hole at $(w_0,m)$ is

\begin{equation}
\frac{2^{\frac{m(m+1)}{2}}}{(m-1)!!}
\begin{vmatrix}
((I+\Delta)^{w_0-1} ((2I+\Delta)^{-1} (x^{j-1})))(1) \\
y_1^{j-1} \\
y_2^{j-1} \\
\dotsb \\
y_{m-1}^{j-1}
\end{vmatrix}
\label{eq:8}
\end{equation}

What would be the weighted number of matchings of a black-edged 
$n\times (m+1)$ Aztec rectangle with dents at $x_1,x_2,\dotsb x_m$
and a black hole at $(w_0,m)$? This is the sum we have to simplify:

\begin{equation}
\sum_{x_i\leq y_i<x_{i+1}} \frac{2^{\frac{m(m+1)}{2}}}{(m-1)!!}
\begin{vmatrix}
((I+\Delta)^{w_0-1} ((2I+\Delta)^{-1} (x^{j-1})))(1) \\
y_1^{j-1} \\
y_2^{j-1} \\
\dotsb \\
y_{m-1}^{j-1}
\end{vmatrix}
\end{equation}

This is equal to $\frac{2^{\frac{m(m+1)}{2}}}{(m-1)!!}$ times

\begin{align}
\sum_{x_i\leq y_i<x_{i+1}} 
&\begin{vmatrix}
((I+\Delta)^{w_0-1} ((2I+\Delta)^{-1} (x^{j-1})))(1) \\
y_1^{j-1} \\
y_2^{j-1} \\
\dotsb \\
y_{m-1}^{j-1}
\end{vmatrix}\\ 
&=
\begin{vmatrix}
((I+\Delta)^{w_0-1} ((2I+\Delta)^{-1} (x^{j-1})))(1) \\
x_1^{j-1}+(x_1+1)^{j-1}+\dotsb + (x_2-1)^{j-1} \\
x_2^{j-1}+(x_2+1)^{j-1}+\dotsb + (x_3-1)^{j-1} \\
\dotsb \\
x_{m-1}^{j-1}+(x_{m-1}+1)^{j-1}+\dotsb + (x_m-1)^{j-1} 
\end{vmatrix} \\
&=  
\begin{vmatrix}
((I+\Delta)^{w_0-1} ((2I+\Delta)^{-1} (x^{j-1})))(1) \\
\frac{1}{j} (B_j(x_2) - B_j(x_1)) \\
\frac{1}{j} (B_j(x_3) - B_j(x_2)) \\
\dotsb \\
\frac{1}{j} (B_j(x_m) - B_j(x_{m-1}))
\end{vmatrix} \\
&=
\begin{vmatrix}
1 & \frac{1}{j-1} B_{j-1}(x_1)\\
0 & ((I+\Delta)^{w_0-1} ((2I+\Delta)^{-1} (x^{j-2})))(1) \\
0 & \frac{1}{j-1} B_{j-1}(x_2) - B_{j-1}(x_1)\\
\dotsb & \dotsb \\
0 & \frac{1}{j-1} B_{j-1}(x_m) - B_{j-1}(x_{m-1})
\end{vmatrix} \\
&=
- \begin{vmatrix}
0 & ((I+\Delta)^{w_0-1} ((2I+\Delta)^{-1} (x^{j-2})))(1) \\
1 & \frac{1}{j-1} B_{j-1}(x_1)\\
0 & \frac{1}{j-1} B_{j-1}(x_2) - B_{j-1}(x_1)\\
\dotsb & \dotsb \\
0 & \frac{1}{j-1} B_{j-1}(x_m) - B_{j-1}(x_{m-1})
\end{vmatrix} \\
&= 
- \begin{vmatrix}
0 & ((I+\Delta)^{w_0-1} ((2I+\Delta)^{-1} (x^{j-2})))(1) \\
1 & \frac{1}{j-1} B_{j-1}(x_1)\\
1 & \frac{1}{j-1} B_{j-1}(x_2)\\
\dotsb & \dotsb \\
1 & \frac{1}{j-1} B_{j-1}(x_m)
\end{vmatrix} 
\end{align}

Since 
$\{((I+\Delta)^{w_0-1} ((2I+\Delta)^{-1} (x^{j-2})))(1)\}_{j=2,3,\dotsb m+1}$
is a linear combination of rows of the form 
$\{k^{j-2}\}_{j=2,3,\dotsb m+1}$, $1\leq k\leq m$, it is enough to show how to
simplify

\begin{equation}
\begin{vmatrix}
0 & k^{j-2} \\
1 & \frac{1}{j-1} B_{j-1}(x_1) \\
1 & \frac{1}{j-1} B_{j-1}(x_2) \\
\dotsb & \dotsb \\
1 & \frac{1}{j-1} B_{j-1}(x_m)
\end{vmatrix}
\end{equation}

This happens to be easier than one would expect. As a particular case of
\[\sum_{n=M}^{N-1} n^q = \frac{B_{q+1}(N)-B_{q+1}(M)}{q+1},\]
we have
\[k^{j-1} = \frac{1}{j} (B_j(k+1)-B_j(k)).\] Therefore

\begin{align}
\begin{vmatrix}
0 & k^{j-2} \\
1 & \frac{1}{j-1} B_{j-1}(x_1) \\
1 & \frac{1}{j-1} B_{j-1}(x_2) \\
\dotsb & \dotsb \\
1 & \frac{1}{j-1} B_{j-1}(x_m)
\end{vmatrix}
&=
\begin{vmatrix}
1-1 & \frac{1}{j-1} (B_{j-1}(k+1)-B_{j-1}(k)) \\
1 & \frac{1}{j-1} B_{j-1}(x_1) \\
1 & \frac{1}{j-1} B_{j-1}(x_2) \\
\dotsb & \dotsb \\
1 & \frac{1}{j-1} B_{j-1}(x_m)
\end{vmatrix}\\
&=
\frac{1}{m!}
\begin{vmatrix}
(k+1)^{j-1}-k^{j-1}\\
x_1^{j-1}\\
x_2^{j-1}\\
\dotsb \\
x_m^{j-1}
\end{vmatrix}\\
&=
\frac{1}{m!}
\begin{vmatrix}
\Delta(x^{j-1})(k)\\
x_1^{j-1}\\
x_2^{j-1}\\
\dotsb \\
x_m^{j-1}
\end{vmatrix}
\end{align}

Thus, we now know that the sequence of elementary column operations
we have to apply to the matrix having $x_i^j$ on its lower rows in order
to make it into a matrix having 
$\{1,\frac{1}{1} B_1(x_i),\frac{1}{2} B_2(x_i),\dotsb \frac{1}{m} B_m(x_i)\}$
on its lower rows transforms the row $\{\Delta(x^{j-1})(k)\}_1^{m+1}$ into the 
row $\{0,k^0,k^1,\dotsb k^m\}$. What row would be transformed into the
row 
\[
\begin{aligned}
\{0,&((I+\Delta)^{w_0-1}((2I+\Delta)^{-1}(x^0))(1),
((I+\Delta)^{w_0-1}((2I+\Delta)^{-1}(x^1))(1),\dotsb ,\\
&((I+\Delta)^{w_0-1}((2I+\Delta)^{-1}(x^{m-1}))(1)\}?
\end{aligned}\]
Let us express $(I+\Delta)^{w_0-1}((2I+\Delta)^{-1})(f)(x)$ 
($f$ any polynomial of degree at most $m-1$) 
in the form $\sum_{k=0}^{m-1} a_k f(x+k)$ for some numbers 
$a_1,a_2,\dotsb a_m$. If rows $y_1,y_2,\dotsb y_r$ are transformed into rows
$z_1,z_2,\dotsb z_r$, respectively, then the row $\sum_{i=1}^r b_i y_i$
must be transformed into the row $\sum_{i=1}^r b_i z_i$. Therefore the
row $\{\sum_{k=0}^{m-1} a_k (\Delta(x^{j-1}))(k+1)\}_1^{m+1}$ is transformed
into 
\[\{ 0,\sum_{k=0}^{m-1} a_k (k+1)^0,\sum_{k=0}^{m-1} a_k (k+1)^1,\dotsb ,
\sum_{k=0}^{m-1} a_k (k+1)^{m-1}\},\]
which is the same as
\[
\begin{aligned}
\{0, &(I+\Delta)^{w_0-1}((2I+\Delta)^{-1}(x^0))(1),
(I+\Delta)^{w_0-1}((2I+\Delta)^{-1}(x^1))(1),\dotsb ,\\
&(I+\Delta)^{w_0-1}((2I+\Delta)^{-1}(x^{m-1}))(1)\}.
\end{aligned}
\]
Now, what is $\sum_{k=0}^{m-1} a_k (\Delta(x^j))(k+1)$ for $0\leq j<m$?
It is 
\[((I+\Delta)^{w_0-1}((2I+\Delta)^{-1}(\Delta(x^j))))(1).\]

Hence the weighted number of matchings of a black-edged $n\times (m+1)$ Aztec 
rectangle with dents at $x_1,\dotsb, x_m$ and a black hole at $(w_0,m)$ is
\begin{equation}
-\frac{2^{\frac{m(m+1)}{2}}}{m!!}
\begin{vmatrix}
((I+\Delta)^{w_0-1}(2I+\Delta)^{-1}\Delta)(x^{j-1})(1)\\
x_1^{j-1}\\
x_2^{j-1}\\
\dotsb \\
x_m^{j-1}
\end{vmatrix}. \label{eq:9}
\end{equation}

In the same way we have arrived at this result, and using it as a base 
case for induction, it is easy to prove the following.

\begin{Lem}\label{Lem:koloro}
The weighted number of matchings of
a white-edged $n\times (m+d_1)$ Aztec rectangle with teeth at
$y_1,\dotsb y_{m+d_1-1}$ and a black hole at $(w_0,m)$ is
\begin{equation}
(-1)^{d_1} \frac{2^{\frac{(m+d_1)(m+d_1+1)}{2} }}{(m+d_1-1)!!}
\begin{vmatrix}
((I+\Delta)^{w_0-1} (2I+\Delta)^{-(d_1+1)} \Delta^{d_1})(x^{j-1})(1) \\
y_1^{j-1} \\
y_2^{j-1} \\
\dotsb \\
y_{m+d_1-1}^{j-1}
\end{vmatrix}
\end{equation}
and that the weighted number of matchings of a black-edged $n\times (m+d_1)$
Aztec 
rectangle with dents at $x_1,x_2,\dotsb ,x_{m+d_1-1}$
 and a black hole at $(w_0,m)$ is
\begin{equation}
(-1)^{d_1} \frac{2^{\frac{(m+d_1)(m+d_1-1)}{2} }}{(m+d_1-1)!!}
\begin{vmatrix}
((I+\Delta)^{w_0-1} (2I+\Delta)^{-d_1} \Delta^{d_1})(x^{j-1})(1) \\
y_1^{j-1} \\
y_2^{j-1} \\
\dotsb \\
y_{m+d_1-1}^{j-1}
\end{vmatrix}.
\end{equation}
\end{Lem}
Notice that we have used the fact that $(2I+\Delta)^{-1}$ and
$\Delta$ commute. They do so because $(2I+\Delta)^{-1}$, on the domain of
polynomials of degree lower than a given bound, is shorthand
for a finite sum of powers of $\Delta$.

\begin{figure}
        \begin{minipage}[b]{0.5\linewidth}
                \centering \includegraphics[height=2in]{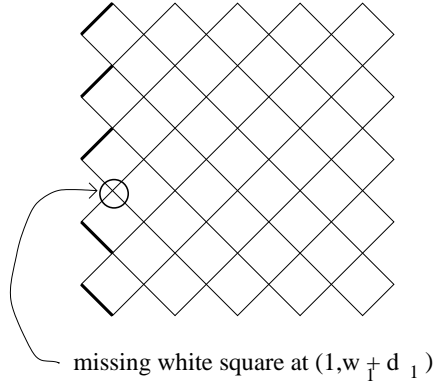}
                \caption{Forced zig-zag pattern}\label{fig:revo}
        \end{minipage}%
\end{figure}

We can now attack our main objective, namely, the enumeration of matchings
of an Aztec diamond with a black hole at $(w_0+d_0,w_1)$ and
a white hole at $(w_0,w_1+d_1)$, where some matchings are counted as
negative matchings. Let us first consider the special case $w_0=1$. 
What happens at white and black vertices $(x,y)$ with $x=1$? The hole
at $(1,w_1+d_1)$ forces a zig-zag pattern covering all those vertices.
(See figure \ref{fig:revo}.) None of the edges covering these vertices
counts towards the weight of the matching. Thus, the weighted number of
matchings of an Aztec diamond of order $n$ with
a black hole at $(1+d_0,w_1)$ and
a white hole at $(1,w_1+d_1)$ is the same as the weighted number of
matchings of a $(n-1)\times n$ white-edged Aztec rectangle with teeth
at $\{1,2,\dotsb n\} $ and a black hole at $(d_0,w_1)$. By
Lemma \ref{Lem:koloro}, this is equal to
\begin{equation}\label{eq:agur}
(-1)^{n-w_1} \frac{2^{\frac{n(n+1)}{2} }}{(n-1)!!}
\begin{vmatrix}
((I+\Delta)^{d_0-1} (2I+\Delta)^{-(n-w_1+1)} \Delta^{n-w_1})(x^{j-1})(1) \\
1^{j-1} \\
2^{j-1} \\
\dotsb \\
(n-1)^{j-1}
\end{vmatrix}.
\end{equation}

The top row is a linear combination of rows of the form
$\{\Delta^k(x^{j-1})(1)\}$, $1\leq j\leq n$. All such rows with 
$k>(n-1)$ vanish. All such rows with $k<(n-1)$ are linear combinations of 
the rows
\[
\{i^{j-1}\}, 1\leq j\leq n, 1\leq i\leq n-1,\label{gaja}
\]
Since an $n$ by $n$ determinant whose last $n-1$ rows are
$\{1^j\},\{2^j\},\dotsb \{(n-1)^j\}$ and whose first row is one of
\ref{gaja} vanishes, we can discard all terms of the form
\[\{\Delta^k(x^{j-1})(1)\}, 1\leq j\leq n, k<(n-1)\]
from the top row of the determinant in
(\ref{eq:agur}). Hence we must consider only
the term of the form 
\[\{C\cdot \Delta^{n-1}(x^{j-1})(1)\} , 1\leq j\leq n\]
in the top row of the determinant. Thus (\ref{eq:agur}) is equal to
\[
\begin{aligned}
\label{eq:malbona}
(-1)^{n-w_1} &\frac{2^{\frac{n(n+1)}{2} }}{(n-1)!!}\cdot 
(C\cdot (-1)^{n-1}\cdot (n-1)!!)\\
&= C\cdot (-1)^{w_1+1} 2^{\frac{n(n+1)}{2}}.
\end{aligned}
\]
times
the coefficient of $\Delta^{n-1}$ in \[(I+\Delta)^{d_0-1}
(2I+\Delta)^{-(n+1-w_1)} \Delta^{n-w_1}.\] Since, in the task of finding 
this coefficient, $\Delta$ plays a purely symbolic role, we might
as well use a symbol which does not denote an operator, such as $x$.
It will also be convenient to use the following notation, due to
Richard Stanley.
\begin{Def} The coefficient of $x^j$ in the formal power series $f(x)$
on $x$ is denoted by $[x^j](f(x))$.
\end{Def}

Our task is to compute 
$[x^{n-1}]((1+x)^{d_0-1} (2+x)^{-(n+1-w_1)} x^{n-w_1}$. We have
\[
\begin{aligned}
\label{eq:seventeen}
[x^{n-1}]((1+x)^{d_0-1} (2+x)^{-(n+1-w_1)} x^{n-w_1}) &=
[x^{w_1-1}]((1+x)^{d_0-1} (2+x)^{-(n+1-w_1)}) \\
&= [x^{w_1-1}]((2+x)^{w_1-1} \frac{(1+x)^{d_0-1}}{(2+x)^n}) \\
&= [(2y)^{w_1-1}]((2+2y)^{w_1-1} \frac{(1+2y)^{d_0-1}}{(2+2y)^n}) \\
&= 2^{-(w_1-1)} \cdot
[y^{w_1-1}]((2+2y)^{w_1-1} \frac{(1+2y)^{d_0-1}}{(2+2y)^n}) \\
&= 2^{-(w_1-1)} \cdot 2^{w_1-1} 2^{-n} \cdot
[y^{w_1-1}]((1+y)^{w_1-1} \frac{(1+2y)^{d_0-1}}{(1+y)^n}) \\
&= 2^{-n} \cdot
[y^{w_1-1}]((1+y)^{w_1-1} \frac{(1+2y)^{d_0-1}}{(1+y)^n})
\end{aligned}
\]

Now, for any formal power series $f(y)$ on $y$,
\begin{equation}
\begin{aligned}
\label{eq:eighteen}
[y^{w_1-1}]((1+y)^{w_1-1}\cdot f(y)) 
&= \sum_{j=0}^{w_1-1} \binom{w_1-1}{j} [y^j](f(y)) \\
&= \sum_{j=0}^{w_1-1} \binom{(w_1-1-j)+j}{j} [y^j](f(y)) \\
&= \sum_{j=0}^{w_1-1} (([z^{w_1-1-j}](\frac{1}{(1-z)^{j+1}})) \cdot
		       ([y^j](f(y)))) \\
&= \sum_{j=0}^{w_1-1} (([z^{w_1-1}](\frac{z^j}{(1-z)^{j+1}})) \cdot
		       ([y^j](f(y)))) \\
&= \sum_{j=0}^{w_1-1} (([z^{w_1}](\frac{z^{j+1}}{(1-z)^{j+1}})) \cdot
		       ([y^{j+1}](y\cdot f(y)))) \\
&= [z^{w_1}]((\frac{z}{1-z})\cdot f(\frac{z}{1-z})) \\
&= [z^{w_1-1}]((\frac{1}{1-z})\cdot f(\frac{z}{1-z})) .
\end{aligned}
\end{equation}

Therefore
\begin{equation}
\begin{aligned}
[x^{w_1-1}]((1+x)^{d_0-1} (2+x)^{-(n+1-w_1)}) &=
2^{-n}\cdot [y^{w_1-1}]((1+y)^{w_1-1} \frac{(1+2y)^{d_0-1}}{(1+y)^n})
\text{(by (\ref{eq:seventeen}))} \\
&= 2^{-n}\cdot 
[z^{w_1-1}](\frac{1}{1-z} \cdot 
\frac{(1+\frac{2z}{1-z})^{d_0-1}}{(1+\frac{z}{1-z})^n}) \\
&= 2^{-n}\cdot 
[z^{w_1-1}](\frac{1}{1-z} \cdot 
(1+\frac{2z}{1-z})^{d_0-1} \cdot (1-z)^n) \\
&= 2^{-n}\cdot 
[z^{w_1-1}](\frac{1}{1-z} \cdot 
(\frac{1+z}{1-z})^{d_0-1} \cdot (1-z)^n) \\
&= 2^{-n}\cdot
[z^{w_1-1}]((1+z)^{d_0-1} \cdot (1-z)^{(n-1)-(d_0-1)})
\label{eq:nineteen}
\end{aligned}
\end{equation}

An expression of this form is called a Krawtchouk polynomial 
(\cite{MS}, p. 130).

We have just proven

\begin{Lem}\label{Lem:Lemur}
The weighted number of matchings of 
an Aztec diamond of order $n$ with a black hole at $(d_0+1,w_1)$ and
a white hole at $(1,w_1+d_1)$, $d_0,d_1\geq 1$, is
\[
(-1)^{w_1+1}\cdot
[z^{w_1-1}]((1+z)^{d_0-1} \cdot (1-z)^{(n-1)-(d_0-1)})\cdot
2^{\frac{n(n-1)}{2}}.
\]
\end{Lem}

We can now work on the general case. Every matching of an
Aztec diamond of order $n$ with a black hole
at $(w_0+d_0,w_1)$ and a white hole at
$(w_0,w_1+d_1)$ can be subdivided
into one, and only one, pair of matchings of the following form.
\begin{enumerate}
\item The first item of the pair is a matching of a white-edged
$n\times (w_1+d_1-1)$ Aztec rectangle with dents at 
$y_1,y_2,\dotsb y_{w_1+d_1-2}$ and with a black hole at $(w_0+d_0,w_1)$.
\item The second item of the pair is a matching of a white-edged
$n\times (n-(w_1+d_1-1))$ Aztec rectangle with dents at 
$z_1,z_2,\dotsb z_{n-(w_1+d_1-1)}$. 
\end{enumerate}
The numbers $y_1,y_2,\dotsb y_{w_1+d_1-2}$,
$z_1,z_2,\dotsb z_{n-(w_1+d_1-1)}$ and $w_0$ are all distinct and cover all
of the interval $\{1,2,\dotsb ,n\}$. Thus every matching can be
described as a partition of $\{1,2,\dotsb ,n\}-\{w_0\}$ into two subsets,
a matching of an Aztec rectangle with the first subset as its set of dents,
and 
a matching of an Aztec rectangle with the second subset as its set of dents.
As we stated in the
previous section, we will weigh each matching $T$
 by a factor of $(-1)^{w(T)}$, where
$w(T)$ is the sum of the number of edges of the form
$((i-1,w_1+1),(i,w_1))$ or $((i,w_1+1),(i,w_1))$ for $1\leq i<w_0+d_0$
and the number of edges of the form
$((i,w_1+d_1),(i,w_1+d_1-1))$ or $((i,w_1+d_1),(i+1,w_1+d_1-1))$
for $1\leq i<w_0$. Thus
the weight $w(T)$ of a matching is equal to the weight of the matching
of the $n\times (w_1+d_1-1)$ white-edged Aztec rectangle which it induces, plus
the number of teeth of this Aztec rectangle whose indices are lower
than $w_0$.

Hence the weighted number of matchings of an Aztec diamond of order $n$ with a
 black hole at $(w_0+d_0,w_1)$ and a white hole at $(w_0,w_1+d_1)$
is equal to the sum of
\begin{equation}
\begin{aligned}
(-1)^{d_1-1} (-1)^{t(w_0,y_1,\dotsb,y_{w_1+d_1-2})}\cdot
\frac{2^{\frac{(w_1+d_1-1)(w_1+d_1)}{2}}}{(w_1+d_1-2)!!}\cdot \\
\begin{vmatrix}
((I+\Delta)^{w_0+d_0-1}((2I+\Delta)^{-d_1}\Delta^{d_1-1}(x^{j-1})))(1) \\
y_1^{j-1} \\
y_2^{j-1} \\
\dotsb \\
y_{w_1+d_1-2}^{j-1} 
\end{vmatrix}
\cdot \\
\frac{2^{\frac{(n-(w_1+d_1-1))((n-(w_1+d_1-1))+1)}{2}}}
	{(n-(w_1+d_1-1)-1)!!} 
\determinant{z_i^{j-1}}_1^{n-(w_1+d_1-1)}
\end{aligned}
\end{equation}
over all partitions of $\{1,2,\dotsb ,n\} - \{ w_0 \}$ into two sets
\[\{y_1,y_2,\dotsb y_{w_1+d_1-2}\}, \{z_1,z_2,\dotsb z_{n-(w_1+d_1-1)}\},\]
where 
$y_1<y_2<\dotsb <y_{w_1+d_1-2}$ and
$z_1<z_2<\dotsb <z_{n-(w_1+d_1-1)}$, and 
$t(w_0,y_1,\dotsb,y_k)$ is equal to how many of $y_1,y_2,\dotsb y_k$ are
less than $w_0$.

We may dispose of the inconvenient $t(w_0,y_1,\dotsb,y_{w_1+d_1-2})$
by expressing the same result as the sum of
\begin{equation}
\begin{aligned}
(-1)^{d_1}
\frac{2^{\frac{(w_1+d_1-1)(w_1+d_1)}{2}}}{(w_1+d_1-2)!!}
\begin{vmatrix}
0 & ((I+\Delta)^{w_0+d_0-1}((2I+\Delta)^{-d_1}\Delta^{d_1-1}(x^{j-2})))(1) \\
\delta_{v_1,w_0} & v_1^{j-2} \\
\delta_{v_2,w_0} & v_2^{j-2} \\
\dotsb \\
\delta_{v_{w_1+d_1-1},w_0} & v_{w_1+d_1-1}^{j-2}
\end{vmatrix}
\cdot \\
\frac{2^{\frac{(n-(w_1+d_1-1))((n-(w_1+d_1-1))+1)}{2}}}
	{(n-(w_1+d_1-1)-1)!!} 
\determinant{z_i^{j-1}}_1^{n-(w_1+d_1-1)}
\end{aligned}
\end{equation}
over all partitions of $\{1,2,\dotsb ,n\}$ into two sets
$\{v_1,\dotsb,v_{w_1+d_1-1}\}$, $z_1,\dotsb z_{n-(w_1+d_1-1)}$,
where $v_1<v_2<\dotsb <v_{w_1+d_1-1}$ and $z_1<z_2<\dotsb<z_{n-(w_1+d_1-1)}$,
and where $\delta_{i,j}$ is equal to $1$ for $i=j$, $0$ for $i\ne j$.

The following lemma follows immediately from Laplace's development of
a determinant (\cite{Turnbull}, pp. 22-25).

\begin{Lem} \label{Lemma:C}
Let us have a determinant $\determinant{b_{i,j}}_1^{m_1+m_2}$, where
$b_{i,j} = c_{i,j}$ for $j\leq m_1$, $b_{i,j} = (-1)^{i-1} d_{i,j-m_1}$ for
$j>m_1$. Then 
\begin{equation}
\determinant{b_{i,j}}_1^{m_1+m_2} =
  (-1)^{\lfloor \frac{m_1+m_2}{2} - \frac{m_1}{2} \rfloor } \cdot
\sum \determinant{c_{x_i,j}}_1^{m_1} \cdot 
     \determinant{d_{y_i,j}}_1^{m_2}
\end{equation}
where the sum is over all partitions of $\{1,2,\dotsb , m_1+m_2\}$
into two sets $\{x_1,x_2,\dotsb , x_{m_1}\}$, $\{y_1,y_2,\dotsb , y_{m_2}\}$,
where $x_1<x_2<\dotsb <x_{m_1}$ and $y_1<y_2<\dotsb < y_{m_2}$,
and where $\lfloor z \rfloor$ is the largest integer less than or
equal to $z$.
\end{Lem}

We can now express our sum over partitions as a determinant. The
weighted number of matchings of an Aztec diamond with two holes is

\begin{equation}
(-1)^{d_1}
\frac{2^{\frac{(w_1+d_1-1)(w_1+d_1)}{2}}}{(w_1+d_1-2)!!}
\frac{2^{\frac{(n-(w_1+d_1-1))((n-(w_1+d_1-1))+1)}{2}}}{(n-(w_1+d_1-1)-1)!!}
(-1)^{\lfloor \frac{n+1}{2} \rfloor - \lfloor \frac{w_1+d_1}{2} \rfloor}
\determinant{d_{i,j}}_1^{n+1},
\label{eq:hoho}
\end{equation} 
where 
\begin{itemize}
\item
$d_{i,1}=0$ for all $1\leq i\leq n+1$, $i\ne w_0+1$,
\item
$d_{w_0+1,1}=1$,
\item
$d_{1,j} = ((I+\Delta)^{w_0+d_0-1}
((2I+\Delta)^{-d_1}\Delta^{d_1-1}(x^{j-2})))(1)$
           for all $1<j\leq w_1+d_1$,
\item
$d_{1,j} = 0$ for all $w_1+d_1<j\leq n+1$,
\item $d_{i,j} = (i-1)^{j-2}$ for $i>1$, $1<j\leq w_1+d_1$,
\item $d_{i,j} = (-1)^{i-1} (i-1)^{j-(w_1+d_1+1)}$ for $i>1$, $j>w_1+d_1$
\end{itemize}

The task ahead is to compute the determinant $\determinant{d_{i,j}}_1^{n+1}$.
For convenience, we will refer to it as $D(w_0,d_0,w_1,d_1)$.
From (\ref{eq:hoho}) and Lemma \ref{Lem:Lemur}, it follows that
\[
\begin{aligned}
D(1,d_0,w_1,d_1) &= 
((-1)^{d_1}
\frac{2^{\frac{(w_1+d_1-1)(w_1+d_1)}{2}}}{(w_1+d_1-2)!!}
\frac{2^{\frac{(n-(w_1+d_1-1))((n-(w_1+d_1-1))+1)}{2}}}{(n-(w_1+d_1-1)-1)!!}
(-1)^{\lfloor \frac{n+1}{2} \rfloor - \lfloor \frac{w_1+d_1}{2} \rfloor})^{-1}
\cdot \\
&(-1)^{w_1+1}\cdot
[z^{w_1-1}]((1+z)^{d_0-1} \cdot (1-z)^{(n-1)-(d_0-1)})\cdot
2^{\frac{n(n-1)}{2}}.
\end{aligned}
\]
We will reduce the general case to the special
case $w_0=0$ by expressing
$D(w_0,d_0,w_1,d_1)-D(w_0-1,d_0,w_1,d_1)$, as the product of
$D(1,d_0+w_0-1,w_1,d_1)$ times something else.

If we take the determinant $\determinant{d_{i,j}}_1^{n+1}$ for
$D(w_0-1,d_0,w_1,d_1)$, and add one to the bases of all powers, we
obtain
\begin{equation}
D(w_0-1,d_0,w_1,d_1) = \determinant{g_{i,j}}_1^{n+1}
\end{equation}
where
\begin{itemize}
\item $g_{i,1} = 0$ for all $1\leq i \leq n+1$, $i\ne w_0$;
\item $g_{w_0,1} = 1$;
\item $g_{1,j} = ((I+\Delta)^{w_0+d_0-1}
((2I+\Delta)^{-d_1}\Delta^{d_1-1}(x^{j-2})))(1)$ for all
$1<j\leq w_1+d_1$
; notice
how we have raised the exponent of $(I+\Delta)$ from $w_0+d_0-2$
to $w_0+d_0-1$;
\item $g_{1,j} = 0$ for all $w_1+d_1<j\leq n+1$;
\item $g_{i,j} = i^{j-2}$ for $i>1$, $1<j\leq w_1+d_1$,
\item $g_{i,j} = (-1)^{i-1} i^{j-(w_1+d_1+1)}$ for $i>1$, $j>w_1+d_1$,
\end{itemize}

The $n$-tuple $\determinant{g_{n+1,j}}_2^{n+1}$ is a linear combination of
the $n$-tuples $\determinant{g_{i,j}}_2^{n+1}$, $2\leq i<n+1$ and of
the $n$-tuple $\{1\}_2^{n+1}$, which is what $g_{1,j}$ would
be if the pattern for $\{g_{n,j}\}_2^{n+1},$
$\{g_{n-1,j}\}_2^{n+1},$\dots $\{g_{2,j}\}_2^{n+1}$ were
continued.

\begin{Lem}\label{Lem:kato}
Let $a_k$ be the
coefficient of $x^k$ in $((x-1)^{w_1+d_1-1} (x+1)^{n-(w_1+d_1)+1})$. Then
\begin{equation}
\sum_{k=0}^{n} a_k (k+1)^{j-2} = 0 \text{ for } 1<j\leq w_1+d_1,
\end{equation}
\begin{equation}
\sum_{k=0}^{n} a_k (-1)^k (k+1)^{j-(w_1+d_1+1)} = 0 
\text{ for } w_1+d_1+1\leq j\leq n+1,
\end{equation}
\end{Lem}
\begin{proof}
\begin{equation}
\Delta^{w_1+d_1-1} x^{j-2} = 0  \text{ for } 1<j\leq w_1+d_1
\end{equation}
implies
\begin{equation}\Delta^{w_1+d_1-1} (\Delta + 2I)^{n-(w_1+d_1)+1} x^{j-2} = 0
\text{ for } 1<j\leq w_1+d_1\end{equation}
If we take the value at $x=1$, we obtain
\begin{equation} \sum_{k=0}^n a_k (k+1)^{j-2} = 0  \text{ for } 1<j\leq w_1+d_1.\end{equation}
Similarly,
\begin{equation}
\Delta^{n-(w_1+d_1)+1} x^{j-(w_1+d_1+1)} = 0 \text{ for } 
w_1+d_1+1\leq j\leq n+1
\end{equation}
implies
\begin{equation}
(\Delta + 2I)^{w_1+d_1-1} (-\Delta)^{n-(w_1+d_1)+1} x^{j-(w_1+d_1+1)} = 0
\text{ for } 
w_1+d_1+1\leq j\leq n+1.
\end{equation}
If we let $H$ be the operator taking $x^j$ to $(x+1)^j$, we can write
\begin{equation}\label{eq:blofeld}
(I+H)^{w_1+d_1-1} (I-H)^{n-(w_1+d_1)+1} x^{j-(w_1+d_1+1)} = 0
\text{ for } 
w_1+d_1+1\leq j\leq n+1
\end{equation}
Clearly the coefficient of
$H^k$ in \[(I+H)^{w_1+d_1-1} (I-H)^{n-(w_1+d_1)+1}\] is equal
to $(-1)^k$ times the coefficient of $H^k$ in 
\[(I-H)^{w_1+d_1-1} (I+H)^{n-(w_1+d_1)+1},\] which is equal to
$(-1)^{w_1+d_1-1}$ times the coefficient of
$H^k$ in 
\[(H-I)^{w_1+d_1-1} (H+I)^{n-(w_1+d_1)+1},\]
that is, $a_k$. Hence
\[ \begin{aligned}
(-1)^{w_1+d_1-1} \sum_{k=0}^n (-1)^k a_k (k+1)^{j-(w_1+d_1+1)} &=
\sum_{k=0}^n ((-1)^{w_1+d_1-1} (-1)^k a_k) (k+1)^{j-(w_1+d_1+1)} \\ &=
(I+H)^{w_1+d_1-1} (I-H)^{n-(w_1+d_1)+1} x^{j-(w_1+d_1+1)}(1) \\ &=
(\Delta + 2I)^{w_1+d_1-1} (-\Delta)^{n-(w_1+d_1)+1} 
  x^{j-(w_1+d_1+1)} \\ &=
(\Delta + 2I)^{w_1+d_1-1} \Delta^{n-(w_1+d_1)+1} 
  x^{j-(w_1+d_1+1)} \\ &= 0
\text{ for } 
w_1+d_1+1\leq j\leq n+1.
\end{aligned}
\]
Therefore
\[
\sum_{k=0}^n a_k (-1)^k (k+1)^{j-(w_1+d_1+1)} = 0 \text{ for } 
w_1+d_1+1\leq j\leq n+1.
\]
\end{proof}

Hence, if we add $a_n$ times the $n$th row, $a_{n-1}$ times the
$(n-1)$th row,\dots $a_2$ times the second row to 
the bottom row of $\determinant{g_{i,j}}_1^{n+1}$, we obtain
the row $\{h_j\}_1^{n+1}$, where 
\begin{itemize}
\item $h_0=[x^{w_0-1}]((x-1)^{w_1+d_1-1} (x+1)^{n+1-(w_1+d_1)})$,
\item $h_j = (-1)\cdot (-1)^{w_1+d_1-1} \text{ for } 1<j\leq w_1+d_1$,
\item $h_j = (-1)\cdot (-1)^{w_1+d_1-1} \text{ for } j> w_1+d_1$.
\end{itemize}

After switching signs and shifting the bottom row of
$D(w_0-1,d_0,w_1,d_1)$ (now  $\{h_j\}_1^{n+1}$) to the second-to-topmost
place, we obtain
\begin{equation}
D(w_0-1,d_0,w_1,d_1) = 
(-1)\cdot (-1)^{w_1+d_1-1}\cdot (-1)^{n-1}\cdot
\determinant{k_{i,j}}_1^{n+1},
\end{equation}
where 
\begin{itemize}
\item $k_{i,1} = 0$ for all $1\leq i\leq n+1$, $i\ne w_0+1$, $i\ne 2$;
\item $k_{2,1} = (-1)\cdot (-1)^{w_1+d_1-1}\cdot
[x^{w_0-1}]((x-1)^{w_1+d_1-1} (x+1)^{n+1-(w_1+d_1)})$;
\item $k_{w_0+1,1} = 1$;
\item $k_{1,j} = 
((I+\Delta)^{w_0+d_0-1}((2I+\Delta)^{-d_1}\Delta^{d_1-1}(x^{j-2})))(1)$
for all $1<j\leq w_1+d_1$;
\item $k_{1,j} = 0$ for $w_1+d_1+1\leq j\leq n+1$;
\item $k_{i,j} = (i-1)^{j-2}$ for $i>1$, $1<j\leq w_1+d_1$;
\item $k_{i,j} = (-1)^i (i-1)^{j-(w_1+d_1+1)}$ for $i>1$, $j\geq w_1+d_1+1$.
\end{itemize}

We multiply the $n-(w_1+d_1)+1$ rightmost columns by $(-1)$, obtaining
\begin{equation}
D(w_0-1,d_0,w_1,d_1) = \determinant{l_{i,j}}_1^{n+1},
\end{equation}
where
\begin{itemize}
\item $l_{i,1} = 0$ for all $1\leq i\leq n+1$, $i\ne w_0+1$, $i\ne 2$;
\item $l_{2,1} = (-1)\cdot (-1)^{w_1+d_1-1}\cdot
[x^{w_0-1}]((x-1)^{w_1+d_1-1} (x+1)^{n+1-(w_1+d_1)})$;
\item $l_{w_0+1,1} = 1$;
\item $l_{1,j} = 
((I+\Delta)^{w_0+d_0-1}((2I+\Delta)^{-d_1}\Delta^{d_1-1}(x^{j-2})))(1)$
for all $1<j\leq w_1+d_1$;
\item $l_{1,j} = 0$ for $w_1+d_1+1\leq j\leq n+1$;
\item $l_{i,j} = (i-1)^{j-2}$ for $i>1$, $1<j\leq w_1+d_1$;
\item $l_{i,j} = (-1)^{i-1} (i-1)^{j-(w_1+d_1+1)}$ for $i>1$, $j\geq w_1+d_1+1$.
\end{itemize}

Therefore
\begin{equation}D(w_0,d_0,w_1,d_1)-D(w_0-1,d_0,w_1,d_1) = 
\determinant{r_{i,j}}_1^{n+1},\end{equation}
where
\begin{itemize}
\item $r_{i,1} = 0$ for all $1\leq i\leq n+1$, $i\ne 2$;
\item $r_{2,1} = (-1)^{w_1+d_1-1}\cdot
[x^{w_0-1}]((x-1)^{w_1+d_1-1} (x+1)^{n+1-(w_1+d_1)})$;
\item $r_{1,j} = 
((I+\Delta)^{w_0+d_0-1}((2I+\Delta)^{-d_1}\Delta^{d_1-1}(x^{j-2})))(1)$
for all $1<j\leq w_1+d_1$;
\item $r_{1,j} = 0$ for $w_1+d_1+1\leq j\leq n+1$;
\item $r_{i,j} = (i-1)^{j-2}$ for $i>1$, $1<j\leq w_1+d_1$;
\item $r_{i,j} = (-1)^{i-1} (i-1)^{j-(w_1+d_1+1)}$ for $i>1$, $j\geq w_1+d_1+1$.
\end{itemize}

This is equal to \begin{equation}
\begin{aligned}
r_{2,1}\cdot D(1,w_0+d_0-1,w_1,d_1) &=
(-1)^{w_1+d_1-1}\cdot 
[x^{w_0-1}]((x-1)^{w_1+d_1-1} (x+1)^{n+1-(w_1+d_1)})\cdot \\
&D(1,w_0+d_0-1,w_1,d_1)\\
&=
[x^{w_0-1}]((1-x)^{w_1+d_1-1}(1+x)^{(n-(w_1+d_1-1))})\cdot \\
&D(1,w_0+d_0-1,w_1,d_1)
\end{aligned}
\end{equation}
Therefore
\begin{equation}
\label{eq:nye}
\begin{aligned}
D(w_0,d_0,w_1,d_1) &=
  (\sum_{j=1}^{w_0-1} D(j+1,d_0,w_1,d_1)-D(j,d_0,w_1,d_1)) \\
  &+ D(1,d_0,w_1,d_1) \\
&=  \sum_{j=1}^{w_0-1} ([x^j](1-x)^{w_1+d_1-1}(1+x)^{n-(w_1+d_1-1)})\cdot
		     D(1,j+d_0,w_1,d_1)\\
&+ D(1,d_0,w_1,d_1)\\
&= \sum_{j=0}^{w_0-1} ([x^j](1-x)^{w_1+d_1-1}(1+x)^{n-(w_1+d_1-1)})\cdot
		     D(1,j+d_0,w_1,d_1)
\end{aligned}
\end{equation}

By (\ref{eq:hoho}), it follows that 
the weighted number of matchings $\sum (-1)^{w(T)}$ of an Aztec diamond of
order n with a black hole at $(w_0+d_0,w_1)$ and a white hole at
$(w_0,w_1+d_1)$ is equal to 
\[
\label{eq:hra}
\begin{aligned}
(-1)^{d_1}
&\frac{2^{\frac{(w_1+d_1-1)(w_1+d_1)}{2}}}{(w_1+d_1-2)!!}
\frac{2^{\frac{(n-(w_1+d_1-1))((n-(w_1+d_1-1))+1)}{2}}}{(n-(w_1+d_1-1)-1)!!}
(-1)^{\lfloor \frac{n+1}{2} \rfloor - \lfloor \frac{w_1+d_1}{2} \rfloor}\cdot \\
&\sum_{j=0}^{w_0-1} ([x^j](1-x)^{w_1+d_1-1}(1+x)^{n-(w_1+d_1-1)})\cdot
		     D(1,j+d_0,w_1,d_1)\\
&= \sum_{j=0}^{w_0-1} (([x^j](1-x)^{w_1+d_1-1}(1+x)^{n-(w_1+d_1-1)})\cdot \\
&((-1)^{d_1}
\frac{2^{\frac{(w_1+d_1-1)(w_1+d_1)}{2}}}{(w_1+d_1-2)!!}
\frac{2^{\frac{(n-(w_1+d_1-1))((n-(w_1+d_1-1))+1)}{2}}}{(n-(w_1+d_1-1)-1)!!}\cdot \\
&(-1)^{\lfloor \frac{n+1}{2} \rfloor - \lfloor \frac{w_1+d_1}{2} \rfloor}
\cdot  D(1,j+d_0,w_1,d_1)) \\
\end{aligned}
\]

The term within parentheses including $D(1,j+d_0,w_1,d_1)$ is equal to
the number of matchings of an Aztec diamond of
order n with a black hole at $(j+d_0+1,w_1)$ and a white hole at
$(1,w_1+d_1)$. By Lemma \ref{Lem:Lemur}, this number is equal to
\[
\begin{aligned}\label{eq:hro}
(-1)^{w_1+1}\cdot 2^{\frac{n(n-1)}{2}}\cdot 
[z^{w_1-1}]((1+z)^{j+d_0-1} \cdot (1-z)^{(n-1)-(j+d_0-1)})
\end{aligned}
\]

Hence $\sum (-1)^{w(T)}$ is equal to
\[
\begin{aligned}
(-1)^{w_1+1}\cdot 2^{-n}\cdot 
\sum_{j=0}^{w_0-1}
(&[x^j]((1-x)^{w_1+d_1-1}(1+x)^{n-(w_1+d_1-1)})\cdot \\
&[z^{w_1-1}]((1+z)^{j+d_0-1} \cdot (1-z)^{(n-1)-(j+d_0-1)})),
\end{aligned}
\]

From this and from (\ref{eq:mandrake}) the result we have sought
follows immediately.

\begin{Prop}\label{prop:soup}
The entry in the inverse of the Kasteleyn matrix of an
Aztec diamond of order $n$ corresponding to a black square at $(w_0+d_0,w_1)$
and a white square at $(w_0,w_1+d_1)$, $d_0,d_1>0$, is
\begin{equation}
\begin{aligned}\label{eq:gp}
(-1)^{d_0+d_1+w_1}\cdot 2^{-n}\cdot 
\sum_{j=0}^{w_0-1}
(&[x^j]((1-x)^{w_1+d_1-1}(1+x)^{n-(w_1+d_1-1)})\cdot \\
&[z^{w_1-1}]((1+z)^{j+d_0-1} \cdot (1-z)^{(n-1)-(j+d_0-1)})),
\end{aligned}
\end{equation}
where $[x^j](p(x))$ is the coefficient of $x^j$ in the polynomial $p(x)$. 
(Alternatively, this can be called the value of the
coupling function of the Aztec
diamond of order $n$ at the black hole $(w_0+d_0,w_1)$ and the white
hole $(w_0,w_1+d_1)$.)
\end{Prop}

In order to deal with the cases $d_0\leq 0$ and $d_1\leq 0$,
we merely need to flip
the Aztec diamond so as to make $d_0,d_1>0$ and
compute the weighted number of tilings in the manner we have described.
Of course, we have to account for the fact that the weighting has to
be computed differently. For $d_0\leq 0$, we also have to express
$D(w_0,d_0,w_1,d_1)$ as
\[
- \sum_{j=w_0}^n D(j+1,d_0,w_1,d_1)-D(j,d_0,w_1,d_1)
\] and not as

\[
D(1,d_0,w_1,d_1) +
\sum_{j=1}^{w_0-1} D(j+1,d_0,w_1,d_1)-D(j,d_0,w_1,d_1)) 
\]
as we did in (\ref{eq:nye}). These are the only two details worth mention
in the otherwise trivial derivation of the following result
from Proposition \ref{prop:soup}.
\begin{Cor}
The coupling function of the Aztec diamond of order
 $n$ at the black square $(w_0+d_0,w_1)$
and the white square $(w_0,w_1+d_1)$ is
\begin{equation}
\begin{aligned}
(-1)^{d_0+d_1+w_1}\cdot 2^{-n}\cdot 
\sum_{j=0}^{w_0-1}
(&[x^j]((1-x)^{w_1+d_1-1}(1+x)^{n-(w_1+d_1-1)})\cdot \\
&[z^{w_1-1}]((1+z)^{j+d_0-1} \cdot (1-z)^{(n-1)-(j+d_0-1)}))
\end{aligned}
\end{equation}
for $d_0>0$,
\begin{equation}
\begin{aligned}
(-1)^{d_0+d_1+w_1}\cdot 2^{-n}\cdot 
(-\sum_{j=w_0}^n
(&[x^j]((1-x)^{w_1+d_1-1}(1+x)^{n-(w_1+d_1-1)})\cdot \\
&[z^{w_1-1}]((1+z)^{j+d_0-1} \cdot (1-z)^{(n-1)-(j+d_0-1)})))
\end{aligned}
\end{equation}
for $d_0\leq 0$.
\end{Cor}
When we take a minor of the inverse Kasteleyn matrix,
the factors $(-1)^{d_0+d_1+w_1}$, multiplied, give the same product
in every term of the expression of the minor in a form such as
$\sum_{\pi} \sgn(\pi ) \prod_{i=1}^k c_{i,\pi(i)}$. Thus we can
leave them out, and our main
result follows.
\begin{Thm}
The probability of a pattern covering white squares $v_1,v_2,\dotsb v_k$
and black squares $w_1,w_2,\dotsb w_k$ of
 an Aztec diamond of order $n$ is equal
to the absolute value of
\[\determinant{c(v_i,w_j)}_{i,j=1,2,\dotsb k}.\]
The {\em coupling function} $c(v,w)$ at white square $v$ and black
square $w$ is 
\[2^{-n} \sum_{j=0}^{x_i-1} \kr(j,n,y_i-1) 
			    \kr({y\prime }_i - 1,n-1,n-(j+{x\prime }_i-x_i))
\]
for ${x\prime }_i > x_i$ and
\[-2^{-n} \sum_{j=x_i}^n \kr(j,n,y_i-1) 
			 \kr({y\prime }_i-1,n-1,n-(j+{x\prime }_i-x_i))
\]
for ${x\prime }_i \leq x_i$, where $(x_i,y_i)$ and $({x\prime}_i,{y\prime}_i)$
are the coordinates of $v$ and $w$, respectively, in the coordinate
system in figure \ref{fig:coor}, and the {\em Krawtchouk polynomial}
$\kr(a,b,c)$ is the coefficient of $x^a$ in $(1-x)^c\cdot (1+x)^{b-c}$.
\end{Thm}
\section{In perspective}

\begin{figure}
\centering 
\rotatebox{-90}{\includegraphics[height=3in]{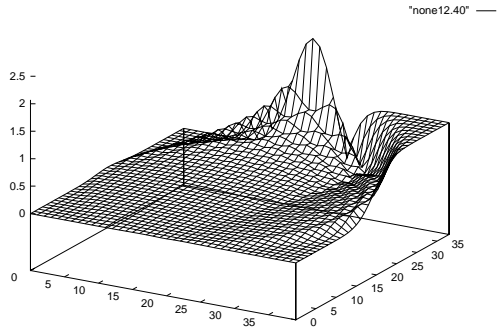}}
\caption{Inverse Kasteleyn matrix entries as a function of $w_0$, $w_1$
for $d_0=1$, $d_1=2$, $n=40$}\label{fig:ha} 
\end{figure}

\begin{figure}
\centering 
\rotatebox{-90}{\includegraphics[height=3in]{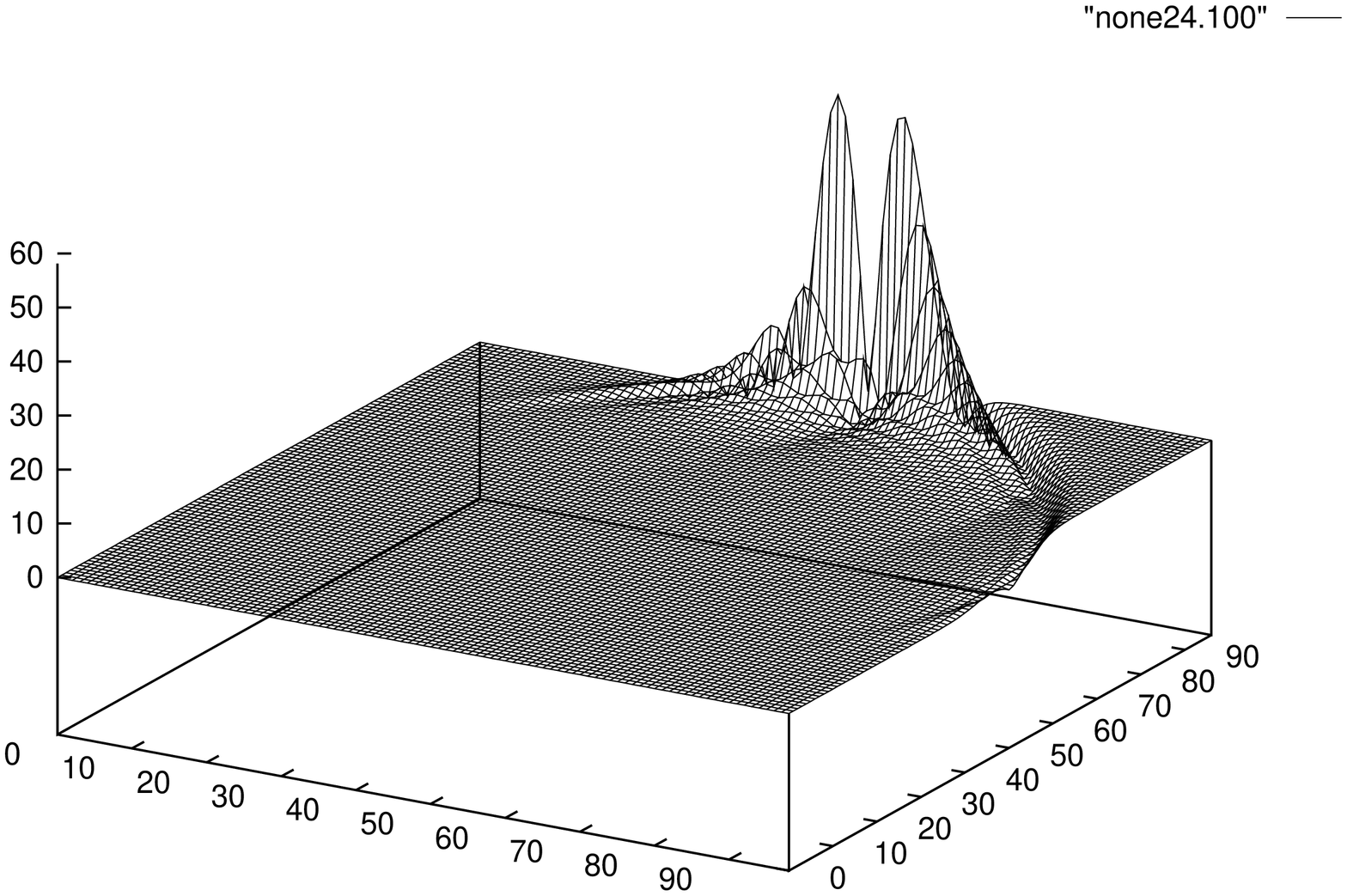}}
\caption{Inverse Kasteleyn matrix entries as a function of $w_0$, $w_1$
for $d_0=2$, $d_1=4$, $n=100$}\label{fig:he}
\end{figure}

Proposition \ref{prop:soup} is valuable in itself, in that it gives us an
efficient algorithm for computing an arbitrary entry in the inverse
Kasteleyn matrix of the Aztec diamond. Figures \ref{fig:ha} to 
\ref{fig:he} show the absolute value of the entry as a function of
$w_0$ and $w_1$, for fixed $d_0$ and $d_1$. As we showed in section 2,
given a pattern consisting of $k$ vertices,
we can compute the probability of its occurence at any point in
a Aztec diamond of given order
by computing $(k/2)^2$ entries of the inverse Kasteleyn
matrix. Thus, for a fixed pattern, the time
required for computing its probability is equal to a constant times
the time required for computing an entry of the inverse Kasteleyn matrix
using Proposition \ref{prop:soup}. Whether computation time grows quadratically
on the order of the Aztec diamond,
or somewhat faster, depends on whether multiplying integers is assumed
to take constant time. What is clear is that we now have an algorithm
that is much more efficient than computing the entries of an inverse
Kasteleyn matrix by actually inverting the matrix or computing minors.

Proposition
\ref{prop:soup} gives us an expression that is more closed than an entry
in the inverse of a Kasteleyn matrix. What do we mean by this?
There are few tools available that would allow us to obtain asymptotic
expressions for a sequence of entries in a sequence of inverses of 
arbitrary matrices. For finding the asymptotics of 
sums such as (\ref{eq:gp}), however, there are
many well-developed analytical techniques.
At the time of this writing, Henry Cohn is working on some minor 
problems involved in applying the saddle-point technique to the asymptotics
of (\ref{eq:gp}). Once he superates these difficulties 
(something that seems to be about to happen), the goals set in the introduction
will have been achieved completely.

\begin{figure}
\centering
\rotatebox{-90}{\includegraphics[height=3in]{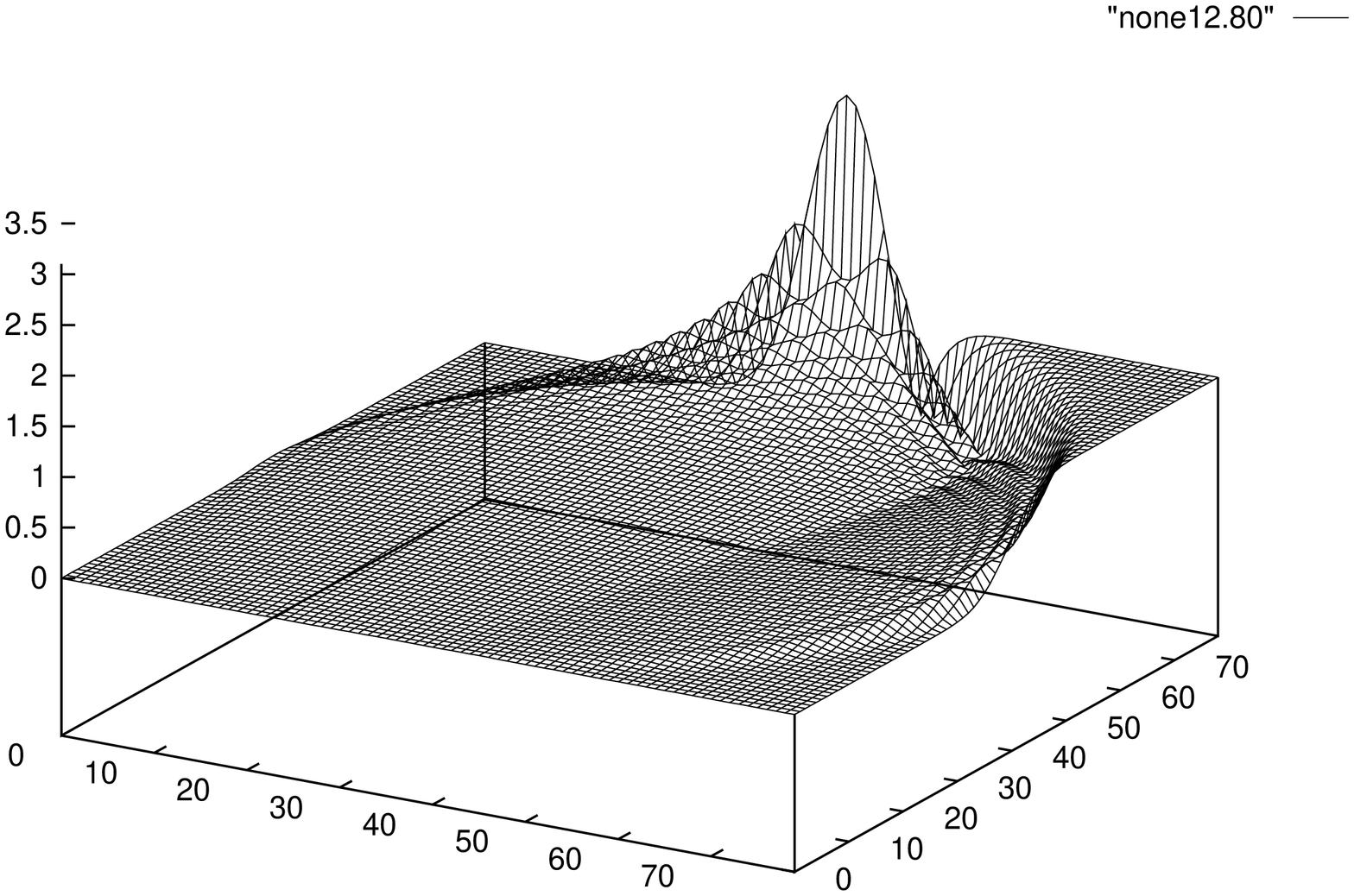}}
\caption{Inverse Kasteleyn matrix entries as a function of $w_0$, $w_1$
for $d_0=1$, $d_1=2$, $n=80$} 
\end{figure}
\begin{figure}
\centering
\rotatebox{-90}{\includegraphics[height=3in]{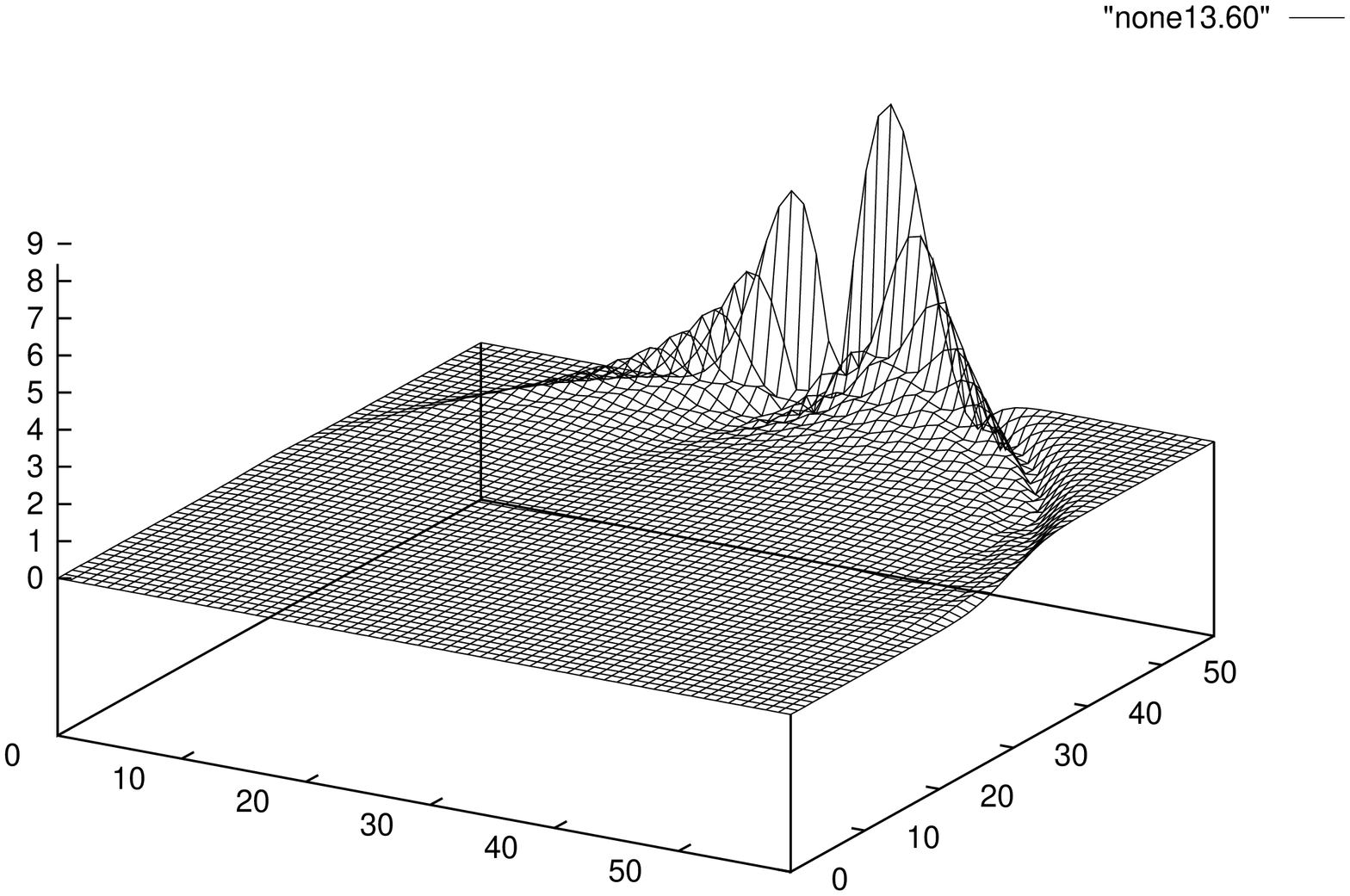}}
\caption{Inverse Kasteleyn matrix entries as a function of $w_0$, $w_1$
for $d_0=1$, $d_1=3$, $n=60$}
\end{figure}
\begin{figure}
\centering
\rotatebox{-90}{\includegraphics[height=3in]{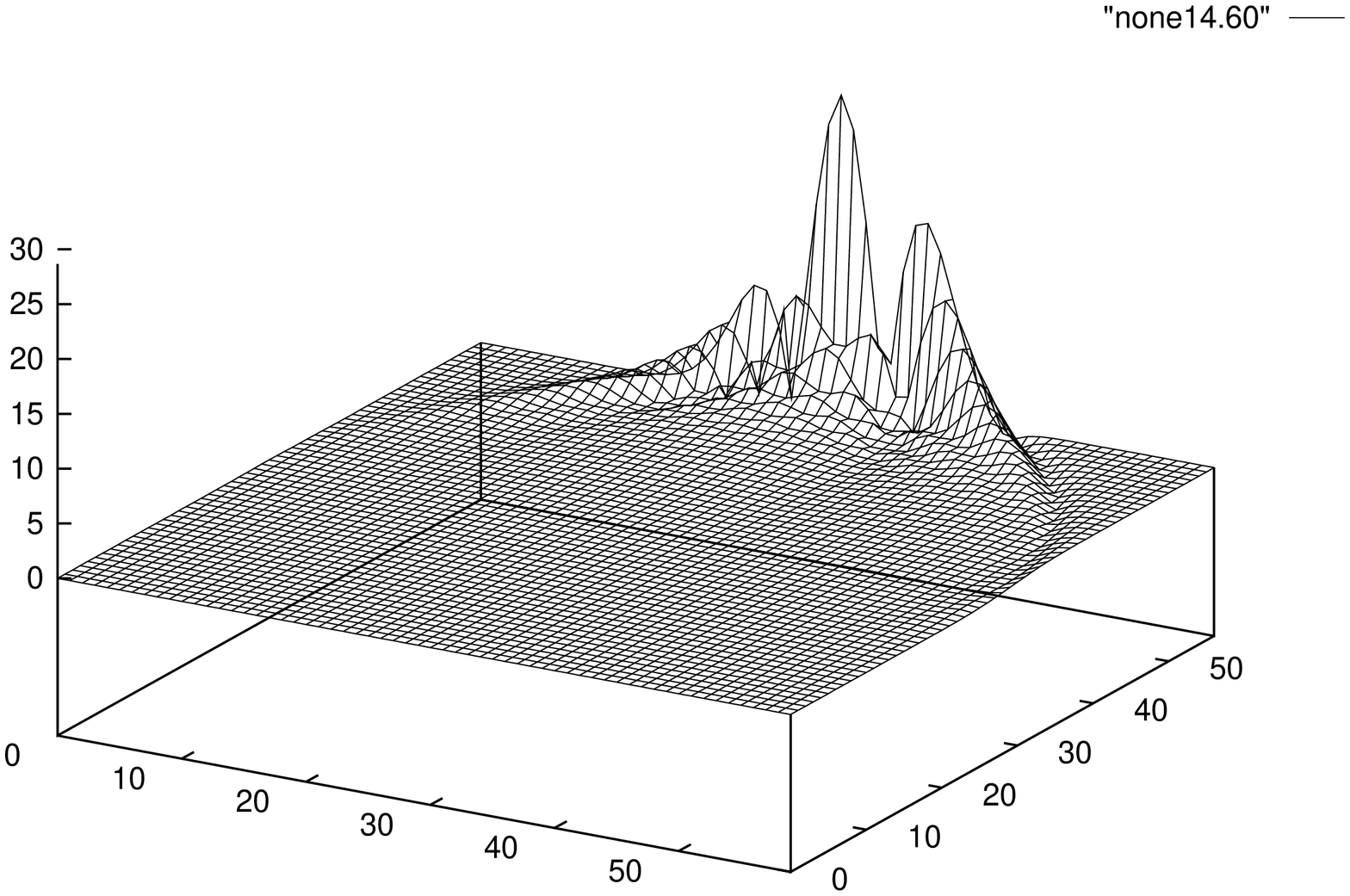}}
\caption{Inverse Kasteleyn matrix entries as a function of $w_0$, $w_1$
for $d_0=1$, $d_1=4$, $n=60$}
\end{figure}
\begin{figure}
\centering
\rotatebox{-90}{\includegraphics[height=3in]{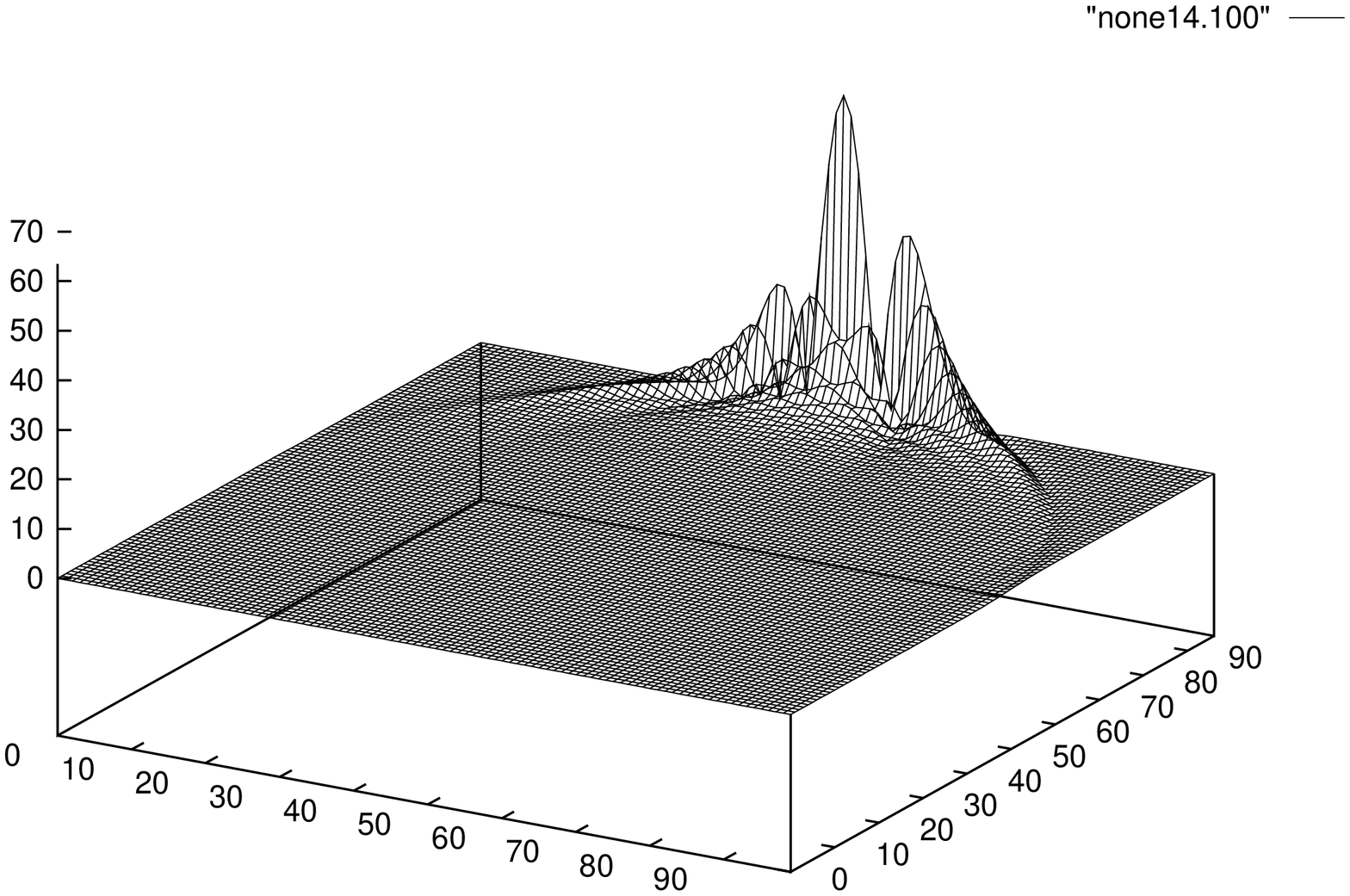}}
\caption{Inverse Kasteleyn matrix entries as a function of $w_0$, $w_1$
for $d_0=1$, $d_1=4$, $n=100$}
\end{figure}
\begin{figure}
\centering
\rotatebox{-90}{\includegraphics[height=3in]{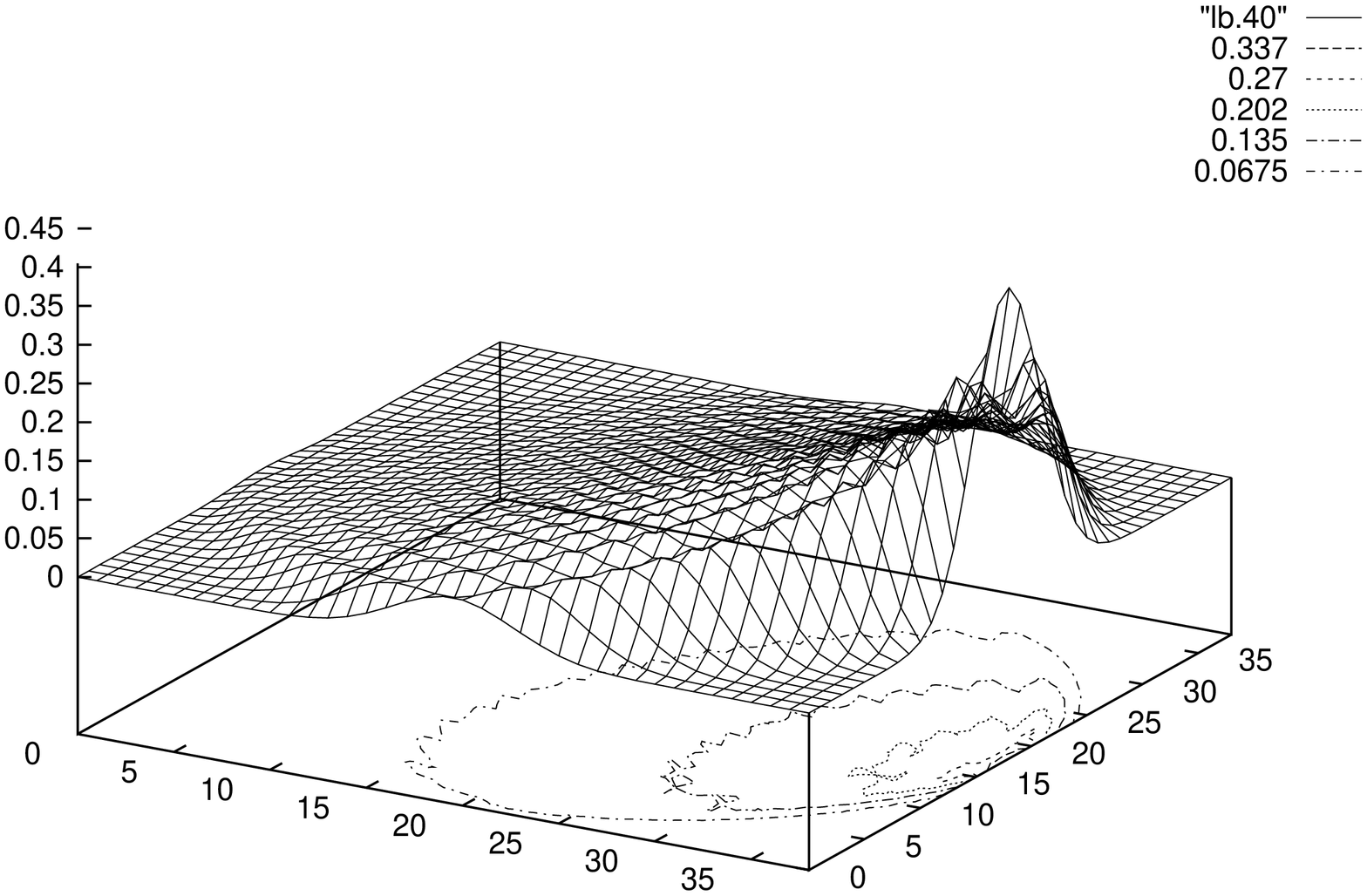}}
\caption{Probability of occurence of the shape in figure \ref{fig:nag}
as a function of position, for the Aztec diamond of side $40$.}
\end{figure}
\begin{figure}
\centering
\rotatebox{-90}{\includegraphics[height=3in]{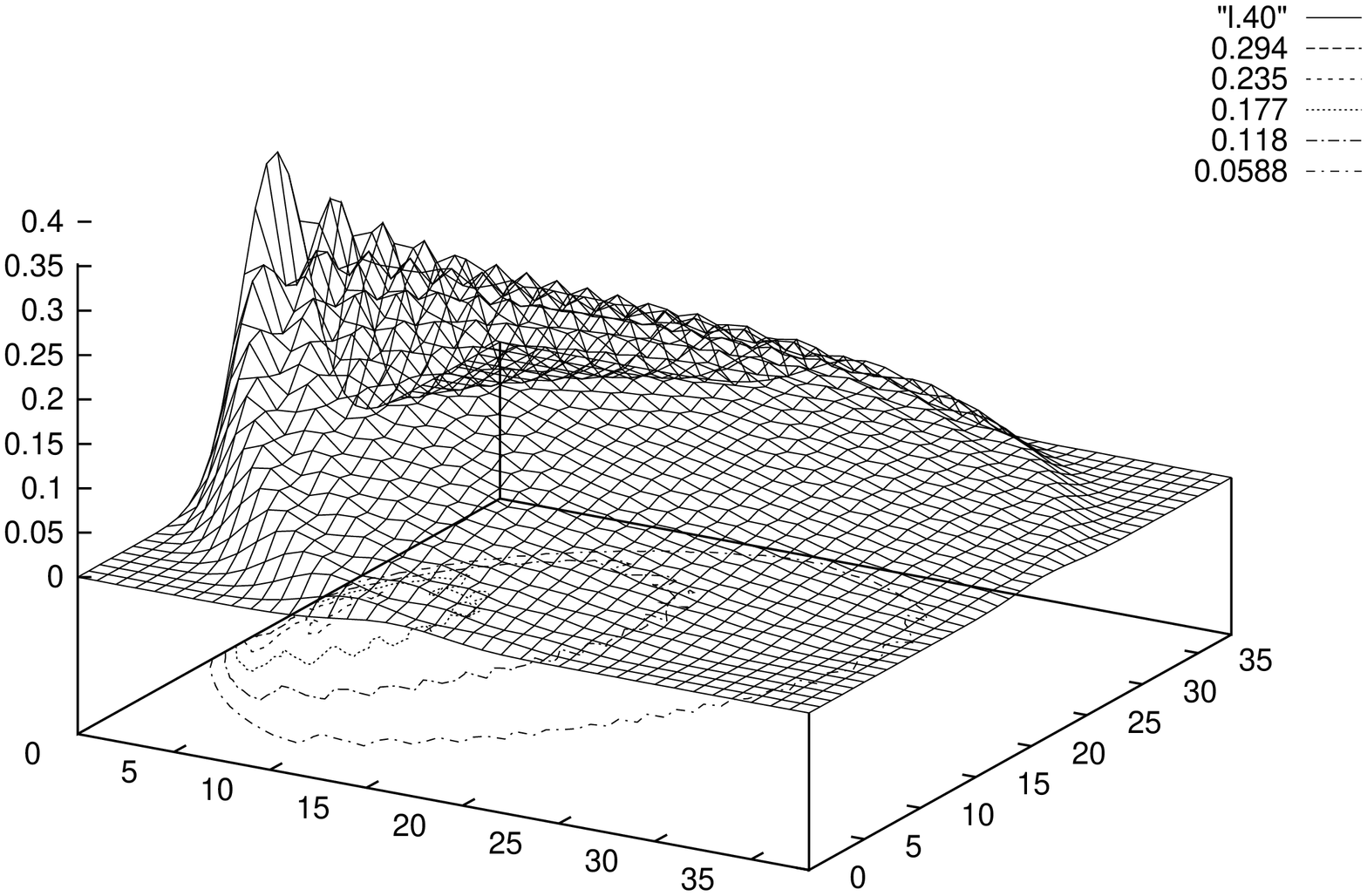}}
\caption{Probability of occurence of the shape in figure \ref{fig:neg}
as a function of position, for the Aztec diamond of side $40$.}
\end{figure}
\begin{figure}
        \begin{minipage}[b]{0.5\linewidth}
                \centering \includegraphics[height=1in]{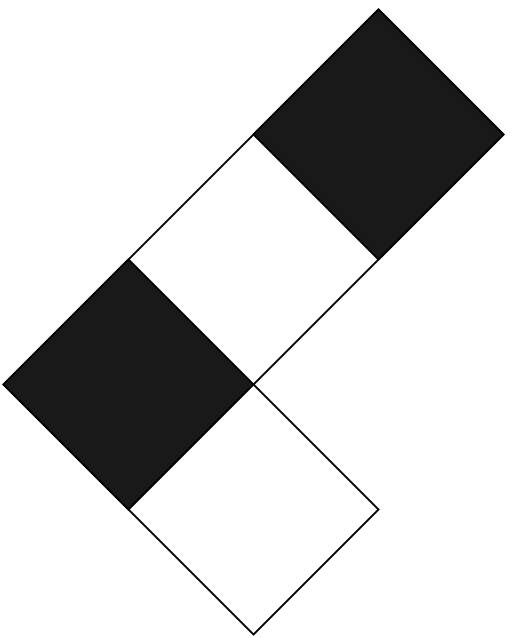}
                \caption{L-shape}\label{fig:nag}
        \end{minipage}%
        \begin{minipage}[b]{0.5\linewidth}
                \includegraphics[height=1in]{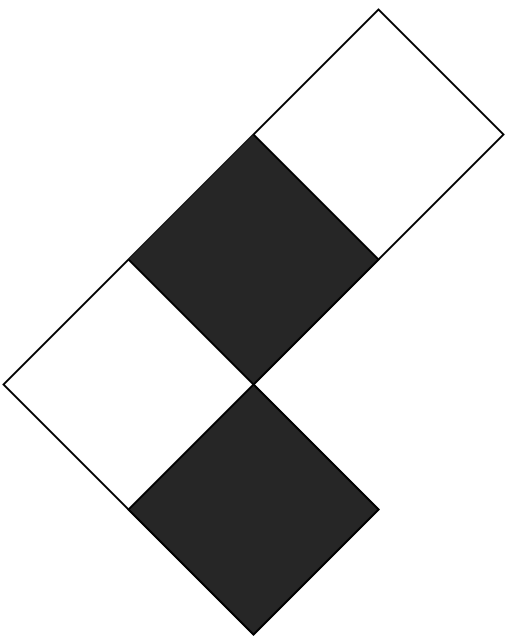}
                \caption{L-shape} \label{fig:neg}   
        \end{minipage}%
\end{figure}
\section{Acknowledgements}
The author would like to thank Matthew Blum for converting the output
of the urban-renewal program \pkg{ren.c} to figures \ref{fig:g} to
\ref{fig:f}. He would also like to thank James Propp and
Ira Gessel  for their support, 
and Henry Cohn, in advance, for finishing his work on the asymptotics of the
herein described results and for writing the still unwritten continuation
to this paper.

\end{document}